\documentclass[11pt, a4paper]{article}

\usepackage{latexsym}
\usepackage{amssymb}
\usepackage{amsbsy}
\usepackage{amsmath}
\usepackage{amsthm}
\usepackage[all]{xy}
\usepackage{mathrsfs}
\usepackage{epsfig}
\usepackage{enumerate}
\usepackage{balance}
\usepackage{arydshln}
\usepackage{graphicx}
\usepackage{graphics}
\usepackage{color}
\usepackage{appendix}
\usepackage{thesis_symbols}

\newcommand{\maphi}{\widehat{\varphi}_T}
\newcommand{\imaphi}{\widehat{\varphi}_{T^{-1}}}

\renewcommand{\ttw}{{q}}
\renewcommand{\ttwe}{{q}}
\renewcommand{\gamma}{r}
\renewcommand{\delta}{k}
\newcommand{\rsp}{{\rm rowspan~}}
\renewcommand{\ann}[1]{{\rm ann}\left(#1\right)}

\title{Novel representation of discrete $n$-D autonomous systems}
\author{Debasattam Pal
\thanks{Department of Electronics and Electrical Engineering,
Indian Institute of Technology Guwahati, Guwahati - 781 039, India.
{\tt debasattam@iitg.ernet.in},
}
\and
Harish K. Pillai
\thanks{Department of Electrical Engineering,
Indian Institute of Technology Bombay, Mumbai - 400 076, India.
{\tt hp@ee.iitb.ac.in}}
}
\date{}

\begin{document}
\maketitle
\begin{abstract}
In this paper we address the problem of
representing solutions of a system of scalar linear partial difference equations akin to
state space equations of $1$-D systems theory.
We first obtain a representation formula for a special class of autonomous systems.
Then we show every autonomous system can be converted into the special ones by a coordinate transformation
on $\Z^n$. Using this conversion we provide representation formula for general autonomous systems.
The representation formula we present can be viewed as multidimensional flow
operators acting on initial conditions. These initial conditions are required to satisfy certain
compatibility conditions. We give a full description of the set of allowable initial conditions.

In our search for a general representation formula, one algebraic result plays a very crucial role.
In this result we show that every quotient ring of the $n$-variable Laurent polynomial ring can be
made a finitely generated faithful module over another Laurent polynomial ring of smaller dimension
by doing a suitable change of coordinates. We call this result a discrete version of Noether's
Normalization Lemma.
\end{abstract}

\section{Introduction}
One of the most crucial discoveries in the theory of dynamical systems is the notion of {\em state variables}.
The state variables serve, simultaneously, two purposes in the analysis and control of dynamical systems.
On the one hand, they relate to the system's memory, or, equivalently, energy storing capability, while on
the other hand, they provide a recursive formula to compute the trajectories in the systems. The hallmark
charateristic of state variables that enables the above-mentioned two benefits is that the representation
of a dynamical system in state variables is a {\em first order} ordinary differential equation
(ODE) of the following form:
\begin{equation}\label{eq:1d:state}
\frac{{\rm  d}x}{{\rm d}t}=f(x,u),~~y = h(x,u).
\end{equation} 
For discrete dynamical systems, that is, systems described by ordinary difference equations (odes), the state equation
takes the form:
\begin{equation}\label{eq:1d:state:dis}
x(k+1)=f(x(k),u(k)),~~y(k) = h(x(k),u(k)).
\end{equation} 
Given the benefits of state variables it is desirable to have systems decribed in the form of equation \eqref{eq:1d:state}.
Unfortunately, often the first principles modelling of dynamical systems end up giving higher order equations \cite{state-rapi}.
However, for a system of ODEs, a higher order mathematical model can always be {\em reduced} to a state variable model (see
\cite{state-rapi} where this was proven for linear time-invariant systems). And, therefore, it is a standard practice to consider
dynamical systems already given by a state variable model.

The situation drastically changes when it comes to systems described by partial differential/difference equations (PDEs/pdes).
For such $n$-D systems\footnote{This way, systems described by ODEs/odes are called $1$-D systems.}, as they are often called,
it is not clear whether a state variable representation similar to equations \eqref{eq:1d:state} or \eqref{eq:1d:state:dis}
always exists. The main hindrance in obtaining such a first order representation for $n$-D systems is that, unlike $1$-D systems,
there is no fixed direction of evolution for $n$-D systems. This problem can be somewhat overcome by letting one of the
independent variables play the role of `time' (see, for example, \cite{curtain,sasane,avelli:roc:rap:10}). This approach naturally
leads to an infinite dimensional
state space \cite{curtain}. A serious drawback of this approach is that there are many $n$-D systems for which such a state space representation
with a special variable treated as time may not even be possible \cite{avelli:roc:rap:10}.

In this paper, we show that for {\em every} discrete $n$-D autonomous system, a first order representation akin to equation
\eqref{eq:1d:state:dis} can be obtained. Our inspiration comes from the special class of autonomous $n$-D systems -- the
{\em strongly autonomous systems}. It is known that a strongly autonomous linear $n$-D system admits a first order
representation of the following type:  
\begin{equation}\label{eq:state:strongly:1}
\left.\begin{array}{lcr}
x(\nu_1+1,\nu_2,\ldots,\nu_n)&=&A_1 x(\nu_1,\nu_2,\ldots,\nu_n)\\
x(\nu_1,\nu_2+1,\ldots,\nu_n)&=&A_2 x(\nu_1,\nu_2,\ldots,\nu_n)\\
&\vdots&\\
x(\nu_1,\nu_2,\ldots,\nu_n+1)&=&A_1 x(\nu_1,\nu_2,\ldots,\nu_n)\\
w(\nu_1,\nu_2,\ldots,\nu_n)&=&Cx(\nu_1,\nu_2,\ldots,\nu_n),
\end{array}\right\}
\end{equation}
where $A_1,A_2,\ldots,A_n$ are real square invertible matrices, $C$ is a real matrix of appropriate size
\cite{roc:wil:89,roc:zam:93}. Note that here,
each state matrix $A_i$ encodes the information about the evolution of the state variable $x$
in the $i^{\rm th}$ direction of the integer grid $\Z^n$, over which the system's trajectories evolve. The state variable $x$
evaluated at any integer $n$-tuple
$(\nu_1,\nu_2,\ldots,\nu_n)\in\Z^n$ is a finite dimensional vector. So the state space here is a finite dimensional vector space.
We show in this paper, that for a class of systems (that strictly contains the class of strongly autonomous systems),
such a first order representation exists. However, the representation has a striking difference from equation
\eqref{eq:state:strongly:1}: here the state space is usually infinite dimensional. We show that for these systems,
there exists a non-negative integer $d<n$ such that the state space is a subspace of the space of all vector valued
discrete trajectories evolving over $\Z^d$ (we call such trajectories $d$-D trajectories). And, the first order representation
provides an $(n-d)$-D evolution on the state space. Interestingly, as we show in this paper, the number $d$ is equal to
a geometric invariant of the autonomous system: it equals the dimension of the {\em characteristic variety} of the system.
As a result of this, the representation of equation \eqref{eq:state:strongly:1} turns out to be a special case, because
the defining property of strongly autonomous systems is that the dimension of their characteristic varieties is zero, therefore,
$d=0$ in this case, and so the state space turns out to be finite dimensional.

We call this special class of autonomous systems {\em strongly relevant of order $d$}. The characterizing property of this
class is a certain algebraic property. The crucial benefit of this algebraic characterization is that every discrete $n$-D
autonomous system can be transformed in to an `equivalent' strongly relevant system by a change of coordinates in
$\Z^n$. This enables us to provide a first order representation for {\em any} discrete autonomous $n$-D system. A crucial role
in this transformation is played by a normalization process, which we call the discrete Noether's normalization lemma
(DNNL), for it is an analogue, applied to Laurent polynomial rings, of the classical Noether's normalization lemma, which is
applied to polynomial rings. 

The above-mentioned first order representations can be utilized to provide explicit solution formulae for a given
discrete autonomous $n$-D system. We provide these formulae as algorithms for solving system of autonomous pdes in this paper.
Interestingly, a parallel can be drawn between these representation formulae and the classical result of Ehrenpreis-Palamodov
integral representation theorem for continuous autonomous $n$-D systems \cite{bjork,obe:90,sturmfels-poly}. The integral
representation formula
can be viewed as an evolutionary equation of initial conditions specified on the characteristic variety (which is a
$d$-D object in $\C^n$). The representation formulae of this paper, too, are evolutionary equations on state spaces,
which are subspaces of $d$-D trajectories. However, a striking difference is that, in our case, since the state
space is just a subspace of $d$-D trajectories, the initial conditions are specified over $\Z^d$ as opposed to the ones
specified over the characteristic variety in the integral representation theorem. Further, our evolutionary equation
does not involve an explicit integration. These two facts make our representation formulae easily implentable on a computer.  

The resuts presented in this paper are generalizations of those in \cite{pal:pillai:sicon}, which deals with the
special case of $n=2$. It is well-known that often such an extension from $n=2$ to a general $n$ is not straight-forward,
owing to the fact that dimension $2$ algebras have certain special features that do not carry over to algebras
with higher dimensions (see \cite{hp:shiva:98} for some of these special features of $2$-variable polynomial rings).
In our case, the main difficulty arises in obtaining the DNNL for higher dimensional rings. Another crucial difference
between $2$-D and $n$-D is that, in case of general $n$, the number $d$ can be anything between $0$ and
$n-1$, while in case of $2$-D the number $d$ can only be $0$ or $1$.

The paper is organized in the following manner: in the next section (Section \ref{sec:prelim}) we discuss the notation
used in this paper, and provide some preliminaries required for the main subject matter of the paper. Then, in Section
\ref{sec:rep:strong}, we consider a special class of discrete autonomous $n$-D systems and show how a first order
representation of such systems can be obtained. Using this representation we then provide an explicit solution formula,
alongwith an algorithm, to compute the trajectories in such systems. In Section \ref{sec:dnnl} we state and prove the discrete
Noether's normalization lemma. Then using this normalization we provide, in Section \ref{sec:rep:gen}, a first order
representation for general autonomous discrete $n$-D systems. As a consequence, we provide an explicit solution formula
for trajectories in general autonomous discrete $n$-D systems, and also an algorithm to compute them. Finally, in Section
\ref{sec:allow}
we provide some interesting insights on the set of allowable initial conditions that are used in the solution formulae.
 
\section{Notation and preliminaries}\label{sec:prelim}
\subsection{Notation}
We use $\R$ and $\C$ to denote the fields of real and complex numbers,
respectively. Consequently, $\R^n$, $\C^n$ denote the $n$-dimensional vector spaces
over $\R$ and $\C$, respectively.
The set of integers is denoted by $\Z$, and $\Z^n$ denotes the set of $n$ tuples of elements in
$\Z$. In this paper, our main object of study is a particular class of multiply indexed sequences of
elements in $\R^\ttw$, for some positive integer $\ttw$.
We denote the set of sequences in $\R$ indexed by $\Z^n$ by the symbol $\W$, {\em i.e.},
$\W:=\{w:\Z^n\rightarrow\R\}$. To denote the set of vector valued sequences indexed by $\Z^n$, we use $\W^\ttw$, {\em i.e.},
$\W^\ttw:=\{w:\Z^n\rightarrow\R^\ttw\}$.
The Laurent polynomial ring in $n$ indeterminates
$\Dx{1},\Dx{2},\ldots,\Dx{n}$, usually written as $\R[\Dx{1}^{\pm 1},\Dx{2}^{\pm 1},\ldots,\Dx{n}^{\pm 1}]$,
will be denoted, in this paper, by $\ash$, and the same in the first $d$
indeterminates $\Dx{1},\Dx{2},\ldots,\Dx{d}$, written as $\R[\Dx{1}^{\pm 1},\Dx{2}^{\pm 1},\ldots,\Dx{d}^{\pm 1}]$,
will be denoted by $\ashd{d}$.
We use $\ash^\ttw$ to denote the free module of rank $\ttw$ over $\ash$, where the elements
of $\ash^\ttw$ are $\ttw$-tuples written as row-vectors. For a set $S$,
we use $S^{m\times n}$ to denote the set of $(m\times
n)$ matrices with entries from the set $S$. We often use ${\rm col}(\bullet)$ to stack up
entries one above the other to make a column vector. The single letter $\uld{}$ is
often used to denote the $n$ tuple $(\Dx{1},\Dx{2},\ldots,\Dx{n})$. Further, for an integer $n$ tuple
$\nu=(\nu_1,\nu_2,\ldots,\nu_n)\in\Z^n$, the symbol $\uld{}^\nu$ denotes the monomial
$\Dx{1}^{\nu_1}\Dx{2}^{\nu_2}\cdots\Dx{n}^{\nu_n}$. In this paper, we follow the bar notation to denote
equivalence classes: for $r(\uld{})\in\ash^\ttw$ and a submodule $\Rmod\subseteq\ash^\ttw$, we use
$\ol{r(\uld{})}$ to denote the equivalence class of $r(\uld{})$ in the quotient module
$\ash^\ttw/\Rmod$.

\subsection{Discrete $n$-D systems}
By discrete $n$-D systems, in this paper, we mean systems described by a set of $n$-D linear partial difference
equations with constant real coefficients. Such partial difference equations are often
described using the $n$-D shift
operators $\Dx{1},\Dx{2},\ldots,\Dx{n}$. These shift operators act on the $n$-D discrete `trajectories' (i.e., real valued
sequences indexed by $\Z^n$) $w~\in~~\R^{\Z^2}$ as follows: for $\nu':=(\nu_1',\nu_2',\ldots.\nu_n'),
\nu:=(\nu_1,\nu_2,\ldots,\nu_n)\in\Z^n$
\begin{equation}\label{eq:shift:action:1}
(\uld{}^{\nu'}w)(\nu)=w(\nu+\nu'),
\end{equation}
where the sum $\nu+\nu'$ is the standard element-wise sum: $\nu+\nu'=(\nu_1+\nu_1',\nu_2+\nu_2',\ldots,
\nu_n+\nu_n')$.
This definition can be extended naturally to define the action of $\ash$, the Laurent polynomial
ring in the shifts, on $\W$. And, likewise, a row-vector
$r(\uld{})=
\left[r_1(\uld{}),r_2(\uld{}), \cdots , r_{\ttw}(\uld{})
\right]$ and a column-vector $w={\rm col}(w_1,w_2,\ldots,w_\ttw)\in\W^\ttw$ we define
$r(\uld{})w:=\sum_{i=1}^{\ttw}r_i(\uld{})w_i$.
This defines the action of the row module $\ash^\ttw$ on
$\W^\ttw$.

In this paper, as a convention, we always consider a vector valued trajectory as a column vector. Consequently, any tuple
of difference operators that
act on a column vector of trajectories and produces scalar trajectories, is written as a row, so that this action
can be compactly written as a `dot' product of a tuple of difference operators and that of trajectories. So, as a
convention, in this paper, elements from a module over an operator ring $\ash$ (or $\ashd{d}$, as the case might be) are
assumed to be row vectors.

\subsection{The kernel representation}
Following Willems' $1$-D terminology \cite{paradigms}, we call the collection of trajectories $w\in\W^\ttw$ that satisfy
a given set of partial difference equations the {\em behavior} of the system, and denote it by $\B$. The above
description of the action of $\ash^\ttw$ on $\W^\ttw$ gives the following representation of  behaviors
of discrete $n$-D systems:
\begin{equation}\label{eq:def:behaviour}
\B:=\{w\in\W^\ttw~|~R(\uld{})w=0\},
\end{equation}
where $R(\uld{})\in\ash^{\bullet\times\ttw}$.
The above equation (\ref{eq:def:behaviour}) is called a {\em kernel
representation}
of $\B$ and written as $\B={\rm ker}(R(\uld{}))$. Note that various different matrices
can have the same kernel. Importantly, all matrices having the same row-span over $\ash$
result in the same behavior. This leads to the following equivalent definition
of behaviors: let $R(\uld{})\in\ash^{\bullet\times\ttw}$ and $\Rmod:=\rsp(R(\uld{}))$,
\begin{equation}\label{eq:def:behaviour:2}
\B(\Rmod):=\{w\in\W^\ttw~|~r(\uld{})w=0 \mbox{ for all }r(\uld{})\in\Rmod\}.
\end{equation}
The submodule $\Rmod$ generated by the rows of a kernel representation matrix is called
{\em the equation module} of $\B$. It was shown in \cite{obe:90} that the submodules of $\ash^\ttw$
and discrete $n$-D behaviors that have $\ttw$ number of dependent variables (that is $w(\nu)\in\R^\ttw$ for all $\nu\in\Z^n$)
are in an inclusion reversing one-to-one correspondence
with each other. We denote the set of all discrete $n$-D behaviors with $\ttw$ number of dependent variables by $\Lw$.

In mathematical modelling of systems, it often turns out to be beneficial to incorporate auxiliary dependent variables called
{\em latent variables} \cite{paradigms}. A representation of the behavior involving both the variables of interest (the
{\em manifest variables} \cite{paradigms}) and the latent variables is called a latent variable representation. The following
equation is one such representation:
\begin{equation}
\B=\left\{w\in\W^\ttw~\vline~\exists~ \ell\in\W^{\ttw'} \mbox{ such that }R(\uld{})w+M(\uld{})\ell=0\right\}.
\end{equation} 
Note that the state variables are a special type of latent variables.

\subsection{Autonomous systems}
In this paper, we deal with a special type of discrete $n$-D systems, namely
{\em autonomous systems}. Among several equivalent definitions of
$n$-D autonomous systems (see \cite{roc:zam:93,hp:shiva:98,val:01}), in this paper, we stick
to the algebraic definition of \cite{hp:shiva:98}, which we state below as
Definition \ref{def:aut}. In Definition \ref{def:aut} and later we need the notion
of {\em characteristic ideal} of a discrete $n$-D system. Let the behavior $\B$ of a discrete
$n$-D system be given by a
kernel representation $\B={\rm ker}(R(\uld{}))$ with $R(\uld{})\in\ash^{\bullet\times\ttw}$.
The {\em characteristic ideal} of $\B$, denoted by $\I(\B)$, is defined as the ideal
of $\ash$ generated by the $(\ttw\times\ttw)$ minors of $R(\uld{})$. If the number of rows in $R(\uld{})$
is less than the number of columns $\ttw$, then $\I(\B)$ is defined to be the zero ideal. 
\BD\label{def:aut}
A discrete $n$-D system is said to be {\em autonomous} if the characteristic ideal $\I(\B)$
is nonzero. Further, an autonomous behavior is said to be {\em strongly autonomous} if
the quotient ring $\ash/\I(\B)$ is a finite dimensional vector space over $\R$.
\ED

\subsection{The quotient module}
Given a behavior $\B={\rm ker}(R(\uld{}))$, let $\Rmod$ be the submodule of $\ash^\ttw$
spanned by the rows of $R(\uld{})$. We define
$$
\M:=\ash^\ttw/\Rmod,
$$
and call it the {\em quotient module} of $\B$. This quotient module $\M$ plays a central role
in this paper. We often let elements from $\M$ act on $\B$. This action is defined as follows:
for $m\in\M$, the action of $m$ on $w\in\B$ is defined to be
the action of a lift of $m$ in $\ash^\ttw$ on $w$. For example, let $r(\uld{})\in\ash^\ttw$ be such that
$\ol{r(\uld{})}=m\in\M$, then\begin{equation}\label{eq:def:q:action}
mw:=r(\uld{})w.
\end{equation}
Note that $m$ may have several distinct lifts in $\ash^\ttw$,
but all of them have the same action on $w\in\B$. This is because two lifts differ by an element in the
equation module $\Rmod$. Therefore, the action of this difference on any $w\in\B$ produces the zero trajectory,
for the action of any element in $\Rmod$ is zero on $\B$. Thus, the definition of action of $\M$ on $\B$ is well-defined.

Now note that it follows from Definition \ref{def:aut}
above that $\B$ is autonomous {\em if and only if} the quotient module $\M$ is a
{\em torsion module}, {\em i.e.}, for every $\ol{m(\uld{})}\in\M$ there exists a $f(\uld{})\in\ash$
such that $\ol{f(\uld{})m(\uld{})}=0\in\M$. In that case we get the following ideal called the
{\em annihilator ideal} of $\M$.
$$
\ann{\M}:=\{f(\uld{})\in\ash~|~f(\uld{})m=0~\mbox{ for all } m\in\M\}.
$$

\section{A special case: $\M$ is a finitely generated module over a smaller ring}\label{sec:rep:strong}
We first consider a special case. Let the $n$-tuple of shift operators be partitioned into two tuples:
$\uld{}=(\uld{1},\uld{2})$, where $\uld{1}=(\Dx{1},\Dx{2},\ldots,\Dx{d})$ and $\uld{2}=(\Dx{d+1},\Dx{d+2},
\ldots,\Dx{n})$. In this section, we consider the case when the quotient module
$\M$ happens to be {\em finitely generated} as a module over the subring $\Ashoneb{}$ (that is, the Laurent
polynomial ring generated by the $d$-tuple of shift operators $\uld{1}$). We denote this smaller ring by $\ashd{d}$.

The motivation for considering this special case comes from the intention to mimic the construction of states
for strongly autonomous systems. Recall that, for strongly autonomous systems, the quotient module $\M$ happens
to be a finite dimensional vector space over $\R$. It is already well-known that, in that case, the $n$-D behavior
admits a state-space representation of the form given by equation \eqref{eq:state:strongly:1} (see \cite{}).
However, there are
many $n$-D autonomous systems that are not strongly autonomous. For such systems the quotient module is an infinite
dimensional vector space over $\R$. Among these systems, there is a special class, for which $\M$ is a finitely
generated module over a smaller ring $\ashd{d}$. We formally define
such systems as {\em strongly relevant system of order $d$}. 

\BD\label{def:str:rel}
Let $\B\in\Lwe$ be an autonomous system. Then $\B$ is said to be {\em strongly relevant of order $d$} if the
the quotient module $\M$ of $\B$ is a finitely generated module over $\ashd{d}$. 
\ED

For example, the following scalar $2$-D system 
$$
\B={\rm ker}\left[\begin{matrix}\Dx{2}^2-2\Dx{2}+1\\\Dx{1}\Dx{2}-\Dx{2}-\Dx{1}+1\end{matrix}\right]
$$
is a strongly relevant system of order $1$. On the other hand, the scalar $2$-D system 
$$
\B={\rm ker}\left[\Dx{1}\Dx{2}-\Dx{2}-\Dx{1}+1\right]
$$
is {\em not} strongly relevant of any order. Note that a strongly relevant system of order $0$ is nothing but a
strongly autonomous system.

The assumption that $\M$ is a finitely generated module over $\ashd{d}$ gurantees the existence
of certain matrices which will be crucially used in obtaining the desired representation formulae. We describe the
construction of these matrices and their properties below. 

\subsection{The matrix of relations, $X(\uld{1})$}\label{sub:sec:1}
The quotient module $\M$ being finitely generated over $\ashd{d}$ implies that we
can find a finite collection of elements from $\M$, say $\{g_1,\ldots,g_\gamma\}$, that
generate $\M$ as an $\ashd{d}$-module. This choice of generators
gives us a surjective $\ashd{d}$-module morphism
\begin{equation}\label{eq:finitely_gen}
\begin{array}{cccc}
\psi:&\ashd{d}^{\gamma}&\twoheadrightarrow&\M\\
 & e_i &\mapsto & {g_i}
\end{array}
\end{equation}
for $1\leqslant i \leqslant \gamma$, where $e_i$ is the $i^{\rm th}$ standard basis
row-vector
$$
\begin{array}{ccc}e_i=& \left[\begin{array}{cccccc}0 & 0 &\cdots & \!\!1 & \cdots & 0\end{array}\right]&\in\ashd{d}^\gamma.\\
 & ~\begin{array}{cccccc}&&&\stackrel{\uparrow}{\mbox{$i^{\rm th}$ position}}&&\end{array}&
\end{array}
$$
Clearly, ${\rm ker}(\psi)$ is a submodule of $\ashd{d}^{\gamma}$. Since
$\ashd{d}^\gamma$ is a Noetherian module, this submodule is finitely generated.
Therefore, we can find a matrix $X(\uld{1})\in\ashd{d}^{\delta\times\gamma}$ whose rows
generate ${\rm ker}(\psi)$ as an $\ashd{d}$-module.
Thus $\M$ turns out to be isomorphic as an $\ashd{d}$-module
with the quotient module $\ashd{d}^\gamma/{\rm rowspan}(X(\uld{1}))$:
\begin{equation}\label{eq:iso1}
\ashd{d}^\gamma/{\rm rowspan}(X(\uld{1}))\simeq\M \mbox{ as $\ashd{d}$-modules}.
\end{equation}
We call $X(\uld{1})$ a {\em matrix of relations} of the generators of $\M$
as an $\ashd{d}$-module.



\subsection{Matrices representing the ${\boldmath (n-d)}$-tuple of shift operators $\Dx{d+1},\ldots,\Dx{n}$}
\label{sub:sec:2}
Now consider the $\ashd{d}$-module morphisms defined for $1\leqslant j \leqslant (n-d)$ as
\begin{eqnarray}
\mu_j:\M &\rightarrow &\M \nonumber\\
{m}&\mapsto&\mu_j\left({m}\right):=
{\Dx{d+j}}\cdot m.\label{eq:multiplication}
\end{eqnarray}
These maps are clearly $\ashd{d}$-module morphisms.
Since $\M$ has been assumed to be finitely generated over $\ashd{d}$,
in order to specify these maps it is enough to specify their actions on
the generators $\{{g_1},\ldots,{g_\gamma}\}$. Suppose that
for $1\leqslant i\leqslant \gamma$ and $1\leqslant j \leqslant (n-d)$
\begin{equation}\label{eq:def:aj1}
\mu_j\left({g_i}\right)={a_{i1,j}(\uld{1})}{g_1}+{a_{i2,j}(\uld{1})}
{g_2}+\cdots + {a_{i\gamma,j}(\uld{1})}{g_\gamma}, 
\end{equation}
where $a_{i1,j}(\uld{1}),a_{i2,j}(\uld{1}),\ldots,a_{i\gamma,j}(\uld{1})\in\ashd{d}$.
We denote by $A_j(\uld{1})$ the following matrix
\begin{equation}\label{eq:def:aj}
A_j(\uld{1})=\left[\begin{array}{cccc}
a_{11,j}(\uld{1})&a_{12,j}(\uld{1})&\cdots & a_{1\gamma,j}(\uld{1})\\
a_{21,j}(\uld{1})&a_{22,j}(\uld{1})&\cdots&a_{2\gamma,j}(\uld{1})\\
\vdots & \vdots & \ddots & \vdots \\ 
a_{\gamma 1,j}(\uld{1})&a_{\gamma 2,j}(\uld{1})&\cdots&a_{\gamma\gamma,j}(\uld{1})
\end{array}\right]\in\ashd{d}^{\gamma\times\gamma}.
\end{equation}
Equation (\ref{eq:def:aj1}) can now be written collectively for all $i$ as
\begin{equation}\label{eq:def:aj2}
\mu_j\left(\left[\begin{array}{c}
{g_1}\\{g_2}\\\vdots\\{g_\gamma}
\end{array}\right]\right)=
{A_j(\uld{1})}\left[\begin{array}{c}
{g_1}\\{g_2}\\\vdots\\{g_\gamma}
\end{array}\right].
\end{equation}
It then follows that multiplication by $\Dx{d+j}$ is {\em represented} by the matrix $A_j(\uld{1})$.
We prove this in the following lemma.

\BL\label{lem:companion}
Suppose $\M$ is a finitely generated module over $\ashd{d}$. Let
$\{{g_1},\ldots,{g_\gamma}\}\subseteq\M$ be a generating set for $\M$
as an $\ashd{d}$-module. Consider $\psi:
\ashd{d}^\gamma\rightarrow\M$ as given in equation
(\ref{eq:finitely_gen}).
Further, for $1\leqslant j\leqslant
(n-d)$, let $\mu_j$ and $A_j(\uld{1})$ be as defined by equations
(\ref{eq:multiplication}) and (\ref{eq:def:aj}), respectively. Then the
following diagram commutes.
\begin{equation}
\begin{array}{rcl}
\ashd{d}^{\gamma} & \stackrel{\psi}{\twoheadrightarrow} & \M\\
A_j(\uld{1})~ \downarrow & & \downarrow~\mu_j\\
\ashd{d}^{\gamma} & \stackrel{\psi}{\twoheadrightarrow} & \M
\end{array}
\end{equation}
where $A_j(\uld{1}):\ashd{d}^{\gamma}\rightarrow\ashd{d}^{\gamma}$ is defined by
right-multiplication by the matrix $A_j(\uld{1})$ to the row-vectors in $\ashd{d}^\gamma$.
\EL
\BP
For any $r(\uld{1})=\left[\begin{array}{cccc}r_1(\uld{1}) & r_2(\uld{1}) & \cdots & r_\gamma(\uld{1})\end{array}\right]\in
\ashd{d}^\gamma$ we want to show $$\mu_j(\psi(r(\uld{1})))=\psi(r(\uld{1})A_j(\uld{1})).$$
From the definition of the $\psi$ map it follows that
\begin{equation}\label{eq:lem:com1}
\psi(r(\uld{1}))=r_1(\uld{1}){g_1}+r_2(\uld{1}){g_2}+\cdots+r_\gamma(\uld{1})
{g_\gamma}.
\end{equation}
Operating $\mu_j$ on both sides of equation (\ref{eq:lem:com1}) and then utilizing equation
(\ref{eq:def:aj2}) we get
\begin{eqnarray*}
\mu_j(\psi(r(\uld{1})))&=&\mu_j\left(r_1(\uld{1}){g_1}+r_2(\uld{1}){g_2}+\cdots+r_\gamma(\uld{1}){g_\gamma}\right)\\
&=& \left[\begin{array}{cccc}{r_1(\uld{1})} & {r_2(\uld{1})} & \cdots & {r_\gamma(\uld{1})}\end{array}\right]
\mu_j\left(\left[\begin{array}{c}
{g_1}\\{g_2}\\\vdots\\{g_\gamma}
\end{array}\right]\right)\\
&=&\left[\begin{array}{cccc}{r_1(\uld{1})} & {r_2(\uld{1})} & \cdots & {r_\gamma(\uld{1})}\end{array}\right]
{A_j(\uld{1})}\left[\begin{array}{c}
{g_1}\\{g_2}\\\vdots\\{g_\gamma}
\end{array}\right]=\psi(r(\uld{1})A_j(\uld{1})).
\end{eqnarray*}
\EP

The matrices  $A_1(\uld{1}),A_2(\uld{1}),\ldots,A_{n-d}(\uld{1})$ are called {\em companion matrices}.
In the lemma below we list out some useful properties of these companion matrices.

\BL\label{lem:prop:f_g}
Consider $\psi$ and
$\mu_j$, $A_j(\uld{1})$, for $1\leqslant j\leqslant (n-d)$, as defined by equations
(\ref{eq:finitely_gen}), (\ref{eq:multiplication}) and (\ref{eq:def:aj}), respectively. 
Then we have the following:
\begin{enumerate}
\item The matrices $A_1(\uld{1}),A_2(\uld{1}),\ldots,A_{n-d}(\uld{1})$ commute pairwise.

\item There exists a set of generators $\{g_1,g_2,\ldots,g_\gamma\}$ for $\M$ as a module over $\ashd{d}$
such that each of the corresponding $A_j(\uld{1})$'s 
is invertible in $\ashd{d}^{r\times r}$ (that is, {\em unimodular} over $\ashd{d}$).

\item With each of $A_j(\uld{1})$'s invertible, we have 
the following commutative diagram for any $(\nu_1,\nu_2,\ldots,\nu_{n-d})\in\Z^{n-d}$,
\begin{equation}\label{eq:commute}
\begin{array}{rcl}
\ashd{d}^{\gamma} & \stackrel{\psi}{\twoheadrightarrow} & \M\\
\left(\prod_{j=1}^{n-d}A_j(\uld{1})^{\nu_j}\right)~ \downarrow & & \downarrow~\left(\mu_1^{\nu_1}\circ\mu_2^{\nu_2}\circ\cdots\circ\mu_{n-d}^{\nu_{n-d}}\right).\\
\ashd{d}^{\gamma} & \stackrel{\psi}{\twoheadrightarrow} & \M
\end{array}
\end{equation}

\item With each of $A_j(\uld{1})$'s invertible, 	
for every $(\nu_1,\nu_2,\ldots,\nu_{n-d})\in\Z^{n-d}$, there exists a matrix
$E(\uld{1})\in\ashd{d}^{\delta \times \delta}$ such that
\begin{equation}\label{eq:cor:companion}
X(\uld{1})\left(\prod_{j=1}^{n-d}A_j(\uld{1})^{\nu_j}\right)=E(\uld{1})X(\uld{1}).
\end{equation}
\end{enumerate}
\EL

Due to its length and technical nature we have given the proof of this lemma in the appendix.

\subsection{The output matrix, $C(\uld{1})$}\label{sub:sec:3}
The output matrix is related with the image of the $\ttw\times\ttw$ identity matrix under the canonical map
$\ash^{\ttw}\twoheadrightarrow\M$. We elaborate on this now. Let $s_i$ denote the $i$th standard basis
(row-)vector for $\ash^{\ttw}$, and $\overline{s_i}$ its image under the canonical map
$\ash^{\ttw}\twoheadrightarrow\M$. Now, consider again the map $\psi:\ashd{d}^\gamma\rightarrow \M$ given by
equation (\ref{eq:finitely_gen}). Since $\psi$ is surjective, for each $s_i\in\ash^\ttw$, with $1\leqslant i
\leqslant \ttw$, there exists $c_i(\uld{1})\in\ashd{d}^\gamma$ such that
$$
\psi(c_i(\uld{1}))=\overline{s_i}.
$$
The {\em output matrix} $C(\uld{1})\in\ashd{d}^{\ttw\times\gamma}$ is obtained by stacking these $c_i(\uld{1})$'s
one below the other. That is,
\begin{equation}\label{eq:def:C}
C(\uld{1}):=\left[\begin{matrix}c_1(\uld{1})\\c_2(\uld{1})\\\vdots\\c_\ttw(\uld{1})\end{matrix}\right]\in
\ashd{d}^{\ttw\times\gamma}.
\end{equation}

\subsection{The first order equation}
The key observation to obtaining a representation formula for strongly relevant systems (that is, $\M$ being
finitely generated over $\ashd{d}$) is the fact that such systems admit a first order representation. This fact is
very much akin to the
well-known notion of state space in $1$D systems theory. However, unlike $1$D systems, the state space
in this case is often {\em not} a finite dimensional vector space. In this subsection, we show the construction
of the state space and the corresponding first order representation of the original system, using the three type of
matrices introduced in Subsections \ref{sub:sec:1}, \ref{sub:sec:2} and \ref{sub:sec:3}.

Recall that for strongly relevant systems, the quotient module $\M$ is a finitely generated $\ashd{d}$-module. This
means we have a matrix of relations $X(\uld{1})$, the companion matrices $A_1(\uld{1}),A_2(\uld{1}),
\ldots,A_{n-d}(\uld{1})$, and the output matrix $C(\uld{1})$. Define the following $d$-D behavior
$$
\X:=\left\{x\in\disww{d}{\gamma}~\vline~X(\uld{1})x=0\right\}.
$$
We call $\X$ the {\em state space}.
 
Now note that, as a consequence of statement 4 of Lemma \ref{lem:prop:f_g}, we have, for all $1\leqslant j\leqslant (n-d)$
\begin{equation}\label{eq:aj:inv}
A_i(\uld{1})\X\subseteq\X.
\end{equation}
That is, $\X$ is invariant under the action of each of the companion matrices $\left\{A_i(\uld{1})\right\}_{i=1,2,\ldots,n-d}$.
Therefore, each $A_i(\uld{1})$ can be thought as an {\em endomorphism} of $\X$. We are going to use these companion
matrices to define a first order $(n-d)$-D system over $\X$. For this purpose, let us denote by $\xvec$ an $(n-d)$-D
trajectory over $\X$. Alternatively, $\xvec$ is a map from $\Z^{n-d}\rightarrow\X$. We now describe
the first order system over $\X$ as given below: for $(\nu_1,\nu_2,\ldots,\nu_{n-d})\in\Z^{n-d}$, 
$$
\begin{array}{ccc}
\xvec(\nu_1+1,\nu_2,\ldots,\nu_{n-d})&=&A_1(\uld{1})\xvec(\nu_1,\nu_2,\ldots,\nu_{n-d})\\
\xvec(\nu_1,\nu_2+1,\ldots,\nu_{n-d})&=&A_2(\uld{1})\xvec(\nu_1,\nu_2,\ldots,\nu_{n-d})\\
&\vdots&\\
\xvec(\nu_1,\nu_2,\ldots,\nu_{n-d}+1)&=&A_{n-d}(\uld{1})\xvec(\nu_1,\nu_2,\ldots,\nu_{n-d}).
\end{array}
$$

\BR\label{rem:id}
Note that $\xvec$ can also be identified as an $n$-D trajectory. This identification is done as follows. Suppose
$\nu=(\nu_1,\nu_2,\ldots,\nu_n)\in\Z^n$ is arbitrary. Then $\xvec$, identified as an $n$-D trajectory,
maps this $\nu$ to $(\xvec(\nu_{d+1},\nu_{d+2},\ldots,\nu_n))(\nu_1,\nu_2,\ldots,\nu_d)$. In what follows
we do this identification implicitly throughout the rest of the paper.
\ER

With the identification mentioned in Remark \ref{rem:id}, it follows that $C(\uld{1})\xvec$ can be viewed
as an element of $\disww{n}{\ttw}$. We summarize these constructions in Table 1 below. 

\begin{center}
{\bf TABLE 1}
\\
\begin{tabular}{|c|}
\hline
$\X:=\left\{x\in\disww{d}{\gamma}~\vline~X(\uld{1})x=0\right\}.$
\\
\hline
$\xvec:\Z^{n-d}\rightarrow\X.$
\\
\hline
For $(\nu_1,\nu_2,\ldots,\nu_{n-d})\in\Z^{n-d}$,\\ 
$\begin{array}{ccc}
\xvec(\nu_1+1,\nu_2,\ldots,\nu_{n-d})&=&A_1(\uld{1})\xvec(\nu_1,\nu_2,\ldots,\nu_{n-d})\\
\xvec(\nu_1,\nu_2+1,\ldots,\nu_{n-d})&=&A_2(\uld{1})\xvec(\nu_1,\nu_2,\ldots,\nu_{n-d})\\
&\vdots&\\
\xvec(\nu_1,\nu_2,\ldots,\nu_{n-d}+1)&=&A_{n-d}(\uld{1})\xvec(\nu_1,\nu_2,\ldots,\nu_{n-d}).
\end{array}$
\\
\hline
$w=C(\uld{1})\xvec.$
\\
\hline
\end{tabular}
\end{center}

With this first order system of equations we are now in a position to state and prove the representation theorem
for the strongly relevant systems of order $d$.

\BT\label{thm:main:spl}
Suppose $\B\in\Lwe$ is a strongly relevant autonomous system of order $d$. That is, the quotient module $\M$ of $\B$
is finitely generated as a module
over $\ashd{d}$. Then there exist a positive integer $\gamma$ and the following matrices:
\begin{itemize}
\item $X(\uld{1})\in\ashd{d}^{\bullet\times\gamma}$,
\item $A_1(\uld{1}),A_2(\uld{1}),\ldots,A_{n-d}(\uld{1})\in\ashd{d}^{\gamma\times\gamma}$,
\item $C(\uld{1})\in\ashd{d}^{\ttw\times\gamma}$,
\end{itemize}
such that $\B$ admits the following representation:
\begin{equation}
\B=\left\{w\in\disww{n}{\ttw}~\vline~\exists ~\xvec\in\disww{n}{\gamma} \mbox{such that equations in Table 1
are satisfied}\right\}
\end{equation}
\ET 

We need the following lemma to prove Theorem \ref{thm:main:spl}. The lemma concerns a test, involving the matrices
$X(\uld{1}),A_1(\uld{1}),\ldots,A_{n-d}(\uld{1}),C(\uld{1})$, as to whether a given vector $r(\uld{})\in\ash^\ttw$
maps to $0$ under the canonical surjection $\ash^\ttw\twoheadrightarrow\M$.  

\BL\label{lem:member}
Let $\M=\ash^{\ttw}/\Rmod$ be a finitely generated module over $\ashd{d}$ with
$\{{g_1},\ldots,{g_\gamma}\}$ as a generating set. Let $\psi:
\ashd{d}^\gamma\rightarrow\M$ be given by equation
(\ref{eq:finitely_gen}) and let $X(\uld{1})\in\ashd{d}^{\delta\times\gamma}$ be a matrix
of relations of $\{{g_1},\ldots,{g_\gamma}\}$.
Further, for $1\leqslant j\leqslant
(n-d)$, let $A_j(\uld{1})$ be as defined by equation
(\ref{eq:def:aj}), and let $C(\uld{1})\in\ashd{d}^{\ttw\times\gamma}$ be as defined in equation \eqref{eq:def:C}.
Suppose $f(\uld{})\in\ash^\ttw$ is given by 
$$
f(\uld{})=\sum_{\nu\in\Z^n}\alpha_\nu\uld{}^\nu
$$
with $\alpha_\nu\in\R^\ttw_{\rm row}$ and only finitely many of $\{\alpha_\nu\}$
nonzero. Then the following statements are equivalent:
\begin{enumerate}
\item $f(\uld{})\in\Rmod$.
\item $\sum_{\nu\in\Z^n}\alpha_\nu\prod_{i=1}^{d}
\Dx{i}^{\nu_i}C(\uld{1})\prod_{j=1}^{n-d}A_j(\uld{1})^{\nu_{d+j}}\in{\rm ker}(\psi)$.
\item There exists a matrix $F(\uld{1})\in\ashd{d}^{\ttw\times\delta}$ such that
$$
\sum_{\nu\in\Z^n}\alpha_\nu\prod_{i=1}^{d}
\Dx{i}^{\nu_i}C(\uld{1})\prod_{j=1}^{n-d}A_j(\uld{1})^{\nu_{d+j}}=F(\uld{1})X(\uld{1}).
$$
\end{enumerate}
\EL
\BP
The proof follows by noting how a typical element, say $f(\uld{})\in\ash^\ttw$, corresponds to a vector in
$\ashd{d}^\gamma$, such that the images of these two elements under the maps $\ash^\ttw\twoheadrightarrow\M$ and
$\psi:\ashd{d}^\gamma\twoheadrightarrow\M$, respectively, coincide with each other. That is, given
$f(\uld{})\in\ash^\ttw$, how we can find out an $r(\uld{1})\in\ashd{d}^\gamma$ such that 
$$
\psi(r(\uld{1}))=\overline{f(\uld{})}.
$$
To get this correspondence, first consider $f(\uld{})$ to be a monomial, that is, $f(\uld{})=\alpha \uld{}^\nu$,
where $\alpha\in\R^{\ttw}_{\rm row}$
and $\nu\in\Z^n$. In this case, the image of $f(\uld{})$ under the canonical map $\ash^\ttw\twoheadrightarrow\M$
is given by
\begin{equation}\label{eq:inter:1}
\overline{f(\uld{})}=\alpha\overline{\uld{}^{\nu}}.
\end{equation}
Let us denote the $i^{\rm th}$ component $\nu$ by $\nu_i$. Then by partitioning $\uld{}$ into $\uld{1}$ and $\uld{2}$
we can rewrite equation \eqref{eq:inter:1} as
\begin{eqnarray}
\overline{f(\uld{})}&=&\overline{\alpha}\overline{\uld{}^{\nu}}\nonumber\\
&=&\overline{\alpha} \prod_{i=1}^{d}\overline{\Dx{i}^{\nu_i}}\prod_{i=d+1}^{n}\overline{\Dx{i}^{\nu_i}}\nonumber\\
&=&\overline{\alpha} \overline{\uld{1}}^{{\boldsymbol \nu}_1}\prod_{i=d+1}^{n}\overline{\Dx{i}}^{\nu_i},\label{eq:inter:2}
\end{eqnarray}
where ${\boldsymbol \nu}_1$ denotes $(\nu_1,\nu_2,\ldots,\nu_d)$. Now note that for $\alpha\in\R^{\ttw}_{\rm row}$,
it follows form the definition of the output matrix $C(\uld{1})$ that
$$
\overline{\alpha}=\psi(\alpha C(\uld{1})).
$$
Therefore, since $\psi$ is an $\ashd{d}$-module morphism, we also get
$$
\overline{\alpha}\overline{\uld{1}}^{{\boldsymbol \nu}_1}=\psi\left(\alpha\uld{1}^{{\boldsymbol \nu}_1} C(\uld{1})\right).
$$
Putting this in equation \eqref{eq:inter:2} we get
\begin{eqnarray}
\overline{f(\uld{})}
&=&\overline{\alpha} \overline{\uld{1}}^{{\boldsymbol \nu}_1}\prod_{i=d+1}^{n}\overline{\Dx{i}}^{\nu_i}\nonumber\\
&=&\psi\left(\alpha\uld{1}^{{\boldsymbol \nu}_1} C(\uld{1})\right)\prod_{i=d+1}^{n}\overline{\Dx{i}}^{\nu_i}\nonumber\\
&=&\left(\mu_1^{\nu_{d+1}}\circ\mu_2^{\nu_{d+2}}\circ\cdots\circ\mu_{n-d}^{\nu_{n}}\right)
\psi\left(\alpha\uld{1}^{{\boldsymbol \nu}_1} C(\uld{1})\right)\nonumber\\
&=&\psi\left(\alpha\uld{1}^{{\boldsymbol \nu}_1} C(\uld{1})\prod_{i=1}^{n-d}{A_i(\uld{1}})^{\nu_{d+i}}\right),\nonumber
\end{eqnarray}
where the last equality follows from statement 3 of Lemma \ref{lem:prop:f_g}.

For the case when $f(\uld{})\in\ash^\ttw$ has more than one terms, we note that $\psi$ is $\ashd{d}$-linear. Therefore,
the above argument can be applied term-wise. Thus, for an $f(\uld{})\in\ash^{\ttw}$ given by 
$$
f(\uld{})=\sum_{\nu\in\Z^n}\alpha_\nu\uld{}^\nu
$$
we get
\begin{equation}
\overline{f(\uld{})}=\psi\left(\sum_{\nu\in\Z^n}\alpha_\nu\prod_{i=1}^{d}
\Dx{i}^{\nu_i}C(\uld{1})\prod_{j=1}^{n-d}A_j(\uld{1})^{\nu_{d+j}}\right)
\end{equation}

The result of the lemma then easily follows by noting that $f(\uld{})\in\Rmod$ if and only if $\overline{f(\uld{})}=0\in
\M$.
\EP

With Lemma \ref{lem:member} we are now in a position to prove Theorem \ref{thm:main:spl}.

\BPref{Theorem \ref{thm:main:spl}}
Suppose $\B\in\Lw$ is an autonomous behavior with quotient module $\M$. Let us define $\B_{\rm aux}$ as the
set of trajectories admitting the representation given by Table 1. That is,
\begin{equation}
\B_{\rm aux}:=\left\{w\in\disww{n}{\ttw}~\vline~\exists ~\xvec\in\disww{n}{\gamma} \mbox{such that equations in Table 1
are satisfied}\right\}.
\end{equation}
Then we have to show the following equality of two sets of trajectories:
$$
\B=\B_{\rm aux}.
$$

{\bf\boldmath (To show $\B\subseteq\B_{\rm aux}$)}
We have to show that $w\in\B$ implies that there exists $\xvec\in\disww{n}{\gamma}$ such that equations in Table 1
are satisfied. For the sake of easy referencing we reproduce the equations given in Table 1: for all
$(\nu_1,\nu_2,\ldots,\nu_{n-d})\in\Z^{n-d}$, $\xvec$ satisfies 
\begin{equation}\label{eq:thm:spl:1}
\xvec(\nu_1,\nu_2,\ldots,\nu_{n-d})\in\X, \mbox{ and}
\end{equation}
\begin{equation}\label{eq:thm:spl:2}
\left.\begin{array}{ccc}
\xvec(\nu_1+1,\nu_2,\ldots,\nu_{n-d})&=&A_1(\uld{1})\xvec(\nu_1,\nu_2,\ldots,\nu_{n-d})\\
\xvec(\nu_1,\nu_2+1,\ldots,\nu_{n-d})&=&A_2(\uld{1})\xvec(\nu_1,\nu_2,\ldots,\nu_{n-d})\\
&\vdots&\\
\xvec(\nu_1,\nu_2,\ldots,\nu_{n-d}+1)&=&A_{n-d}(\uld{1})\xvec(\nu_1,\nu_2,\ldots,\nu_{n-d}),\end{array}\right\}
\end{equation}
with 
\begin{equation}\label{eq:thm:spl:3}
w=C(\uld{1})\xvec.
\end{equation}
In order to get this $\xvec$ suppose $\{g_1,g_2,\ldots,g_\gamma\}\subseteq\M$ is a generating set for $\M$
as a module over $\ashd{d}$. Now, let us define for $w\in\B$ the following $\gamma$-tuple of $n$-D trajectories:
\begin{equation}\label{eq:def:x}
\xvec:=\left[\begin{array}{c}g_1 \\g_2 \\\vdots\\g_\gamma \end{array}\right]w.
\end{equation}
We claim that the $\xvec$ defined above satisfies equations \eqref{eq:thm:spl:1}, \eqref{eq:thm:spl:2} and
\eqref{eq:thm:spl:3}.

We first show that $\xvec(\nu_1,\nu_2,\ldots,\nu_{n-d})\in\X$ for all $(\nu_1,\nu_2,\ldots,\nu_{n-d})\in\Z^{n-d}$.
Note that it follows from the defining equation \eqref{eq:def:x} of $\xvec$ that
$$
\xvec(\nu_1,\nu_2,\ldots,\nu_{n-d})=\left[\begin{array}{c}g_1 \\g_2 \\\vdots\\g_\gamma \end{array}\right]w(\bullet,\nu_1,
\nu_2,\ldots,\nu_{n-d}).
$$
We make the matrix of relations $X(\uld{1})$ act on $\xvec(\nu_1,\nu_2,\ldots,\nu_{n-d})$ to get
\begin{eqnarray*}
X(\uld{1})\xvec(\nu_1,\nu_2,\ldots,\nu_{n-d})&=&X(\uld{1})\left[\begin{array}{c}g_1 \\g_2 \\\vdots\\g_\gamma \end{array}\right]
w(\bullet,\nu_1,\nu_2,\ldots,\nu_{n-d})\\
&=&\overline{X(\uld{1})}\left[\begin{array}{c}g_1 \\g_2 \\\vdots\\g_\gamma \end{array}\right]
w(\bullet,\nu_1,\nu_2,\ldots,\nu_{n-d})\\
&=&\psi\left(X(\uld{1})\right)w(\bullet,\nu_1,\nu_2,\ldots,\nu_{n-d})\\
&=& 0,
\end{eqnarray*}
because ${\rm rowspan}(X(\uld{1}))={\rm ker}(\psi)$. Therefore, $\xvec$ satisfies equation \eqref{eq:thm:spl:1}.

We now show that $\xvec$ satisfies the system of equations \eqref{eq:thm:spl:2}. Let $i\in\{1,2,\ldots,(n-d)\}$
and $(\nu_1,\nu_2,\ldots,\nu_{n-d})\in\Z^{n-d}$ be arbitrary. Recall how $\xvec$ is treated as an element in
$\disww{n}{\gamma}$; see Remark \ref{rem:id}. Under this identification, it follows that
\begin{equation}\label{eq:kijani}
\xvec(\nu_1,\nu_2,\ldots,\nu_i+1,\ldots,\nu_{n-d})=\Dx{d+i}\xvec(\nu_1,\nu_2,\ldots,\nu_{n-d}).
\end{equation}
Applying the defining equation \eqref{eq:def:x} of $\xvec$ to equation \eqref{eq:kijani} above we get
\begin{eqnarray*}
\xvec(\nu_1,\nu_2,\ldots,\nu_i+1,\ldots,\nu_{n-d})&=&\Dx{d+i}\xvec(\nu_1,\nu_2,\ldots,\nu_{n-d})\\
&=&\Dx{d+i}\left[\begin{array}{c}g_1 \\g_2 \\\vdots\\g_\gamma \end{array}\right]
w(\bullet,\nu_1,\nu_2,\ldots,\nu_{n-d})\\
&=&\overline{\Dx{d+i}}\left[\begin{array}{c}g_1 \\g_2 \\\vdots\\g_\gamma \end{array}\right]
w(\bullet,\nu_1,\nu_2,\ldots,\nu_{n-d})\\
&=&\mu_i\left(\left[\begin{array}{c}g_1 \\g_2 \\\vdots\\g_\gamma \end{array}\right]\right)
w(\bullet,\nu_1,\nu_2,\ldots,\nu_{n-d})\\
&=&A_i(\uld{1})\left[\begin{array}{c}g_1 \\g_2 \\\vdots\\g_\gamma \end{array}\right]
w(\bullet,\nu_1,\nu_2,\ldots,\nu_{n-d}),
\end{eqnarray*}
where the last equation follows from equation \eqref{eq:def:aj2}. Therefore, we get 
$$
\xvec(\nu_1,\nu_2,\ldots,\nu_i+1,\ldots,\nu_{n-d})=A_i(\uld{1})\xvec(\nu_1,\nu_2,\ldots,\nu_i,\ldots,\nu_{n-d})
$$
for all $1\leqslant i \leqslant (n-s)$ and $(\nu_1,\nu_2,\ldots,\nu_{n-d})\in\Z^{n-d}$. That is,
$\xvec$ satisfies equation \eqref{eq:thm:spl:2}.

Finally, we have to show that $\xvec$ and $w$ satisfy equation \eqref{eq:thm:spl:3}. This follows from the definition
of the output matrix $C(\uld{1})$. Multiplying both sides of equation \eqref{eq:def:x} by $C(\uld{1})$ we get
\begin{eqnarray*}
C(\uld{1})\xvec &=& C(\uld{1})\left[\begin{array}{c}g_1 \\g_2 \\\vdots\\g_\gamma \end{array}\right]w\\
&=& \psi\left(C(\uld{1})\right)w = \overline{I_{\ttw}}w=w.
\end{eqnarray*}
This concludes the $\B\subseteq\B_{\rm aux}$.

{\bf \boldmath (To show $\B\supseteq\B_{\rm aux}$)}
We now show that if we obtain a $w$ by solving the system of equations in Table 1, that is, equations
\eqref{eq:thm:spl:1}, \eqref{eq:thm:spl:2}, \eqref{eq:thm:spl:3}, then that $w$ is in $\B$. Let
$\Rmod\subseteq\ash^{\ttw}$ be the equation module of $\B$. Then, it is enough to show that the action of every element
$r(\uld{})\in\Rmod$ on any $w\in\B_{\rm aux}$ produces the zero trajectory.

Now, let $r(\uld{})$ be any element in $\Rmod$ given in an expanded form as
$$
r(\uld{})=\sum_{\nu\in\Z^n}\alpha_\nu\uld{}^{\nu},
$$
where $\alpha_\nu\in\R^{\ttw}_{\rm row}$ and all of them but finitely many are zero. Once again, we partition the shift
operators
into two, namely $\{\Dx{1},\ldots,\Dx{d}\}$ and $\{\Dx{d+1},\ldots,\Dx{n}\}$, to rewrite $r(\uld{})$ as
$$
r(\uld{})=\sum_{\nu\in\Z^n}\alpha_\nu\prod_{i=1}^{d}\Dx{i}^{\nu_i}\prod_{i=1}^{n-d}\Dx{d+i}^{\nu_{d+i}},
$$
where $(\nu_1,\nu_2,\ldots,\nu_n)=\nu$. Then, this $r(\uld{})$ is made to act upon an arbitrary $w\in\B_{\rm aux}$ to
get the following:
\begin{eqnarray}
r(\uld{})w &=& \left(\sum_{\nu\in\Z^n}\alpha_\nu\prod_{i=1}^{d}\Dx{i}^{\nu_i}\prod_{i=1}^{n-d}\Dx{d+i}^{\nu_{d+i}}\right)w
\nonumber\\
&=&\left(\sum_{\nu\in\Z^n}\alpha_\nu\prod_{i=1}^{d}\Dx{i}^{\nu_i}\prod_{i=1}^{n-d}\Dx{d+i}^{\nu_{d+i}}\right)C(\uld{1})
\xvec\nonumber\\
&=& \left(\sum_{\nu\in\Z^n}\alpha_\nu\prod_{i=1}^{d}\Dx{i}^{\nu_i}C(\uld{1})\prod_{i=1}^{n-d}\Dx{d+i}^{\nu_{d+i}}\right)
\xvec\nonumber\\
&=& \left(\sum_{\nu\in\Z^n}\alpha_\nu\prod_{i=1}^{d}\Dx{i}^{\nu_i}C(\uld{1})
\left(\prod_{i=1}^{n-d}\Dx{d+i}^{\nu_{d+i}} \xvec\right)\right).\label{eq:thm:spl:4}
\end{eqnarray}
Before proceeding, note that 
$$
\Dx{d+i}^{\nu_{d+i}}\xvec=A_i(\uld{1})^{\nu_{d+i}}\xvec,
$$
because $\xvec$ satisfies the first order system equations \eqref{eq:thm:spl:2}. Putting this in equation
\eqref{eq:thm:spl:4} above we get
\begin{eqnarray}
r(\uld{})w &=& \left(\sum_{\nu\in\Z^n}\alpha_\nu\prod_{i=1}^{d}\Dx{i}^{\nu_i}C(\uld{1})
\left(\prod_{i=1}^{n-d}A_i(\uld{1})^{\nu_{d+i}} \xvec\right)\right)\nonumber\\
&=& \left(\sum_{\nu\in\Z^n}\alpha_\nu\prod_{i=1}^{d}\Dx{i}^{\nu_i}C(\uld{1})
\prod_{i=1}^{n-d}A_i(\uld{1})^{\nu_{d+i}} \right)\xvec.\label{eq:thm:spl:5}
\end{eqnarray}
Since $r(\uld{})\in\Rmod$, by Lemma \ref{lem:member}, there exists $F(\uld{1})\in\ashd{d}^{\ttw\times k}$ such that
$$
\sum_{\nu\in\Z^n}\alpha_\nu\prod_{i=1}^{d}\Dx{i}^{\nu_i}C(\uld{1})
\prod_{i=1}^{n-d}A_i(\uld{1})^{\nu_{d+i}} = F(\uld{1})X(\uld{1}).
$$
Making use of this observation in equation \eqref{eq:thm:spl:5} above we get that
\begin{eqnarray*}
r(\uld{})w &=&\left(\sum_{\nu\in\Z^n}\alpha_\nu\prod_{i=1}^{d}\Dx{i}^{\nu_i}C(\uld{1})
\prod_{i=1}^{n-d}A_i(\uld{1})^{\nu_{d+i}} \right)\xvec\\
&=& F(\uld{1})X(\uld{1})\xvec=0,
\end{eqnarray*}
where the last equation follows from the fact that $\xvec$ satisfies equation \eqref{eq:thm:spl:1}, that is,
$$
\xvec(\nu_1,\ldots,\nu_{n-d})\in\X={\rm ker}(X(\uld{1}))
$$ for all $(\nu_1,\ldots,\nu_{n-d})\in\Z^{n-d}$.
Therefore, $w\in\B$, and this concludes the $\B\supseteq\B_{\rm aux}$ part of the proof.
\EP

We illustrate this result in the following example.
\BEx\label{ex:spl:1}
Let us consider the $3$-D scalar system given by the following kernel representation:
$$
\B={\rm ker~}\left[\begin{matrix}\Dx{3}^2-2\Dx{3}+1\\
                                 \Dx{2}^2-2\Dx{2}+1\\
                                 \Dx{1}\Dx{3}-\Dx{1}-\Dx{2}-\Dx{3}+2\end{matrix}\right].
$$
In this case, clearly, $\ash=\R[\Dx{1},\Dx{1}^{-1},\Dx{2},\Dx{2}^{-1},\Dx{3},\Dx{3}^{-1}]$. The role of
the equation module $\Rmod$ is played by the ideal $\a=\langle \Dx{3}^2-2\Dx{3}+1,\Dx{2}^2-2\Dx{2}+1,
\Dx{1}\Dx{3}-\Dx{1}-\Dx{2}-\Dx{3}+2\rangle$. It can be checked that the quotient ring $\M=\ash/\a$ is
finitely generated as a module over $\ashd{1}=\R[\Dx{1},\Dx{1}^{-1}]$. Therefore, $\B$ is strongly relevant
of order $1$. Let us compute the required matrices. As in the main text, let $A_1(\Dx{1}),A_2(\Dx{2})$
represent multiplications by $\Dx{1},\Dx{2}$, respectively, in $\M$. First, we need to fix a generating set for
$\ash/\a$ as a module over $\ashd{1}$. The set $\{\overline{1},\overline{\Dx{2}},\overline{\Dx{3}},
\overline{\Dx{2}\Dx{3}}\}$ meets the requirement for an appropriate generating set. The matrices $A_1(\Dx{1})$
and $A_2(\Dx{1})$, in this case, are given by
$$
\begin{array}{cc}
A_1(\Dx{1})=\left[\begin{matrix}0 & 1 & 0 & 0\\-1 & 2 & 0 & 0\\
                                0 & 0 & 0 & 1\\0 & 0 & -1 & 2 \end{matrix}\right],&
A_2(\Dx{1})=\left[\begin{matrix}0 & 0 & 1 & 0\\0 & 0 & 0 & 1\\
                                -1 & 0 & 2 & 0\\0 & -1 & 0 & 2 \end{matrix}\right].
\end{array}
$$
Note that both $A_1(\Dx{1})$ and $A_2(\Dx{1})$ are invertible over $\ashd{1}$. A matrix of relations, $X(\Dx{1})$,
for the chosen generating set is given by
$$
X(\Dx{1})=\left[\begin{matrix}1 & -1 & -1 & 1\\-\Dx{1}+2 & -1 & \Dx{1}-1 & 0\end{matrix}\right].
$$
It is easily checked that statements 3 and 4 of Lemma \ref{lem:prop:f_g} are satisfied here, that is,
$$
\begin{array}{cc}
X(\Dx{1})A_1(\Dx{1})=\left[\begin{matrix}1 & 0\\\Dx{1}-1 & 1\end{matrix}\right]X(\Dx{1}), &
 X(\Dx{1})A_1(\Dx{1})=\left[\begin{matrix}1 & 0\\-1 & 1\end{matrix}\right]X(\Dx{1}). 
\end{array}
$$  
The state space $\X$ is thus given by
$$
\X={\rm ker~}X(\Dx{1})\subseteq \disww{}{4}.
$$
The output matrix $C(\Dx{1})$ turns out to be
$$
C(\Dx{1})=\left[\begin{matrix}1 & 0 & 0 & 0\end{matrix}\right].
$$
Now, according to Theorem \ref{thm:main:spl}, $w\in\B$ if and only if there exists $\xvec:\Z^{2}\rightarrow\X$ such that
\begin{equation}\label{eq:janine:1}
\left.\begin{array}{rcl}
\xvec(\nu_1+1,\nu_2)&=&A_1(\Dx{1})\xvec(\nu_1,\nu_2)\\
\xvec(\nu_1,\nu_2+1)&=&A_2(\Dx{1})\xvec(\nu_1,\nu_2),
\end{array}\right\}
\end{equation}
and $w=C(\Dx{1})\xvec$.
\EEx

\BR
It is interesting to note that in Example \ref{ex:spl:1}, by utilizing the first order representation, an explicit
solution formula for $w\in\B$ can be given. Note that every $\xvec:\Z^2\rightarrow\X$ satisfying equation
\eqref{eq:janine:1} can be obtained as
$$
\xvec(\nu_1,\nu_2)=A_1(\Dx{1})^{\nu_1}A_2(\Dx{1})^{\nu_2}x,
$$
where $x\in\disww{}{4}$ satisfies $X(\Dx{1})x=0$, that is, $x\in\X$. Now note that $w=C(\Dx{1})\xvec$. Therefore,
by the identifiaction as per Remark \ref{rem:id},
$$
w(\nu_1,\nu_2,\nu_3)=(C(\Dx{1})\xvec(\nu_2,\nu_3))(\nu_1).
$$
With this, now we can make use
of the explicit solution formula for $\xvec$ to obtain a similar solution formula for $w$. This is given by
\begin{equation}\label{eq:for:1}
w(\nu_1,\nu_2,\nu_3)=\left[\begin{matrix}1 & 0 & 0 & 0\end{matrix}\right]
\left[\begin{matrix}0 & 1 & 0 & 0\\-1 & 2 & 0 & 0\\
                                0 & 0 & 0 & 1\\0 & 0 & -1 & 2 \end{matrix}\right]^{\nu_2}
\left[\begin{matrix}0 & 0 & 1 & 0\\0 & 0 & 0 & 1\\
                                -1 & 0 & 2 & 0\\0 & -1 & 0 & 2 \end{matrix}\right]^{\nu_3}x(\nu_1)
\end{equation}
We utilize this solution formula of equation \eqref{eq:for:1} in Algorithm \ref{spl-algo} below.
\ER

As a corollary to Theorem \ref{thm:main:spl} we get an equivalent latent-variable representation for every
strongly relevant autonomous system of order $d$.

\BC\label{cor:main:spl}
Let $\B\in\Lwe$ be an autonomous system that is strongly relevant of order $d$. Then there exist
\begin{itemize}
\item a positive integer $\gamma$,
\item a matrix $X(\uld{1})\in\ashd{d}^{\bullet\times\gamma}$,
\item a set of $n-d$ square matrices $A_1(\uld{1}),\ldots,A_{n-d}(\uld{1})\in\ashd{d}^{\gamma\times\gamma}$,
\item a matrix $C(\uld{1})\in\ashd{d}^{\ttw\times\gamma}$
\end{itemize}
such that $\B$ admits the following latent variable representation:
\begin{equation}\label{eq:cor:lat}
\B=\left\{w\in\disww{n}{\ttw}~\vline~\exists ~\xvec\in\disww{n}{\gamma} \mbox{ satisfying }\\
\left[\begin{array}{c:c}X(\uld{1}) & 0 \\\hdashline
\Dx{d+1}I_{\gamma}-A_1(\uld{1}) & 0\\
\Dx{d+2}I_{\gamma}-A_2(\uld{1}) & 0\\
\vdots & \vdots \\
\Dx{n}I_{\gamma}-A_{n-d}(\uld{1}) & 0\\\hdashline
-C(\uld{1}) & I_{\ttw}
\end{array}\right]\left[\begin{matrix}\xvec \\ w\end{matrix}\right]=0
\right\}
\end{equation}
\EC

\subsection{An algorithm for obtaining all solutions}
Theorem \ref{thm:main:spl} provides a useful algorithm for obtaining all solutions to a strongly relevant autonomous
system.
We describe this algorithm now. The algorithm involves two steps, which we call {\em regularization} and
{\em recursion}.

Given a strongly relevant autonomous system of order $d$, we first obtain the matrices
$X(\uld{1})$, $A_1(\uld{1}),\ldots,A_{n-d}(\uld{1})$, and $C(\uld{1})$. We call this process the {\em regularization}
step. This step requires finding out a suitable set of generators for the finitely generated $\ashd{d}$-module $\M$.
One has to be careful in choosing the generating set so as to ensure that the companion matrices $A_1(\uld{1}),\ldots,
A_{n-d}(\uld{1})$ are all invertible over $\ashd{d}$. This can be implemented on a computer by standard computer
algebra packages. In this paper, we do not go into the details of issues related to such an implementation.

The regularization step leads to the {\em recursion} step, where the companion matrices are used recursively on
trajectories in $\X={\rm ker}(X(\uld{1}))$ to obtain solutions to the original system. The recursion indeed gives the
solutions because of the first order system equations listed in Table 1. 

In order to see that the recursion gives an explicit solution formula for trajectories in the originial system, first
note that solutions to the first order system defined on $\X$ by equations given in Table 1 can be written explicitly
as follows:
\begin{equation}\label{eq:chhata:1}
\xvec(\nu_1,\nu_2,\ldots,\nu_{n-d})=\prod_{i=1}^{n-d}A_i(\uld{1})^{\nu_i}x,
\end{equation}
for arbitrary $(\nu_1,\nu_2,\ldots,\nu_{n-d})\in\Z^{n-d}$, where $x\in\disww{d}{\gamma}$ satisfies $X(\uld{1})x=0$.
Since $w=C(\uld{1})\xvec$, recalling the identification of $\xvec$ as a trajectory in $\disww{n}{\ttw}$ we get
$$
w(\nu_1,\nu_2,\ldots,\nu_n)=\left(C(\uld{1})\xvec(\nu_{d+1},\nu_{d+2},\ldots,\nu_n)\right)(\nu_1,\nu_2,\ldots,\nu_d),
$$
for $(\nu_1,\nu_2,\ldots,\nu_n)\in\Z^n$. Utilizing equation \eqref{eq:chhata:1} now we get
\begin{equation}\label{eq:chhata:2}
w(\nu_1,\nu_2,\ldots,\nu_n)=\left(C(\uld{1})\prod_{i=1}^{n-d}A_i(\uld{1})^{\nu_{d+i}}x\right)(\nu_1,\nu_2,\ldots,\nu_d).
\end{equation}
Algorithm \ref{spl-algo} below elaborates how this equation \eqref{eq:chhata:2} can be used to obtain {\em all}
solutions to a given strongly relevant autonomous system. 

\BA{\bf (Solution formula for strongly relevant systems of order $d$)}\label{spl-algo}

\begin{tabular}{|l|}
\hline
{\bf Regularization}~~\vline\\
\hline
{\bf Input:} An autonomous system $\B={\rm ker}(R(\uld{}))$, where $R(\uld{})\in\ash^{\bullet\times\ttw}$.\\
\hline
{\bf Computation:}\\ 
Obtain a generating set $\mathcal{G}:=\{g_1,\ldots,g_{\gamma}\}\subseteq\M$.\\
Obtain the matrix of relations $X(\uld{1})\in\ashd{d}^{k\times\gamma}$.\\
Obtain the companion matrices $A_1(\uld{1}),\ldots,A_{n-d}(\uld{1})\in\ashd{d}^{\gamma\times\gamma}$\\
(ensure that these matrices are invertible in $\ashd{d}$).\\
Obtain the output matrix $C(\uld{1})\in\ashd{d}^{\ttw\times\gamma}$.\\
\hline
{\bf Output:} The matrices $X(\uld{1})$, $A_1(\uld{1}),\ldots,A_{n-d}(\uld{1})$, and $C(\uld{1})$.\\
\hline
{\bf Construction of $\X$}~~\vline\\
\hline
{\bf Input:} The matrix $X(\uld{1})$.\\
\hline
{\bf Computation:} Solve $X(\uld{1})x=0$ for $x\in\disww{d}{\gamma}$.\\
\hline
{\bf Output:} The set $\X=\{x\in\disww{d}{\gamma}~\vline~ X(\uld{1})x=0\}$.\\
\hline
{\bf Recursion}~~\vline\\
\hline
{\bf Input:} An $x\in\X$, and the matrices $A_1(\uld{1}),\ldots,A_{n-d}(\uld{1})$, and $C(\uld{1})$.\\
\hline
{\bf Computation:}${}$\\
{\bf for} $\nu=(\nu_1,\nu_2,\ldots,\nu_n)\in\Z^n$ set\\
~~~~~~~~~~~~~~~~~~~~
$
w(\nu)=\left(C(\uld{1})\prod_{i=1}^{n-d}A_i(\uld{1})^{\nu_{d+i}}x\right)(\nu_1,\nu_2,\ldots,\nu_d).
$
\\
\hline
{\bf Output:} The sequence $w:=\{w(\nu)\}_{\nu\in\Z^n}$.\\
\hline
\end{tabular}
\EA

\BR
From the algorithm above, it makes sense to call the $d$-D trajectories $x\in\X$ as {\em initial conditions}.
Indeed, the companion matrices act on $x$ as a $(n-d)$-D {\em flow} operator to produce all the solutions in the
original $n$-D system $\B$. We call this condition that $x\in\X$ that is $X(\uld{1})x=0$ as {\em compatibility
condition}. Note that ideally the compatibility condition must be satisfied at every step of the recursion (see
Theorem \ref{thm:main:spl}). However, it turns out that checking for compatibility only once -- at the beginning
of the recursion -- is enough. This is because the set of trajectories satisfying the compatibility conditions,
that is, $\X={\rm ker}(X(\uld{1}))$ is {\em invariant} under each of the companion matrices $A_1(\uld{1}),\ldots,
A_{n-d}(\uld{1})$ (as a consequence of Lemma \ref{lem:prop:f_g}; see equation \eqref{eq:aj:inv}). This is perhaps
the most crucial benefit of the regularization step that it makes the compatibility condition and the recursion
independent of each other, plus it ensures that once the compatibility condition is satisfied, it is automatically
satisfied in every subsequent step of the recursion.
\ER

\section{Discrete version of Noether's normalization}\label{sec:dnnl}
There are many autonomous $n$-D systems that are not strongly relevant. For example, the scalar system
$\B={\rm ker}(\Dx{1}\Dx{2}-\Dx{1}-\Dx{2}+1)$. Here the quotient module is not a finitely generated
module over either $\R[\Dx{1},\Dx{1}^{-1}]$ or $\R[\Dx{2},\Dx{2}^{-1}]$; neither is it a finite dimensional
vector space over $\R$. Therefore, this $\B$ is not a strongly relevant system of any order. In fact, autonomous
systems are more often {\em not} strongly relevant. For such systems, clearly, the method of solution presented above
becomes inapplicable. However, we shall see in this paper that by a normalization process every autonomous system
can be rendered to an equivalent strongly relevant system. We call this normalization process the {\em discrete
Noether's normalization lemma (DNNL)}; for it is an analogous version (for Laurent polynomial rings) to the well-known
Noether's normalization lemma (for polynomial rings). The result was discovered by E. Noether (1882--1935),
who first presented it in her 1926 paper,
entitled ``Der Endlichkeitsatz der Invarianten endlicher linearer Gruppen der Charakteristik $p$", and
used it for studying invariant theory of finite groups over fields of arbitrary characteristics. For detailed analysis
and applications of this result see \cite{eisenbud}.

Like the classical Noether's normalization, our DNNL, too, crucially involves a coordinate transformation. However,
the space on which we need to do the coordinate transformation differs because we are dealing with Laurent polynomial
rings. It turns out that for DNNL we need to do a coordinate transformation on the indexing set $\Z^n$.


\subsection{Coordinate transformations in $\Z^n$ and their effects on $\ash$ and $\ash^{\ttw}$}
The idea of change of coordinates and its effects on a behavior are not new (see \cite{zam:98,
val:01}), however, we use this idea to a completely new end.
By a coordinate change we mean a $\Z$-linear map from $\Z^n$ to
itself of the form
$$
\begin{array}{rcc}
T:\Z^n & \rightarrow & \Z^n\\
{\rm col}(\nu_1,\nu_2,\ldots,\nu_n)=:\nu & \mapsto & T\nu,
\end{array}
$$
where $T\in\Z^{n\times n}$ is a matrix whose determinant is $\pm 1$. Such matrices are called {\em $\Z$-unimodular}.
Note that ${\rm det}(T)=\pm 1$ ensures that the columns of $T$ span the whole of
$\Z^n$ as a $\Z$-module. Such a coordinate transformation $T$ induces the
following $\R$-algebra homomorphism from $\ash$ to itself.
\begin{equation}\label{eq:def:phi}
\begin{array}{rcc}
\varphi_{T}:\ash &\rightarrow& \ash\\
\uld{}^\nu &\mapsto& \uld{}^{T\nu}
\end{array}
\end{equation}
for all $\nu\in\Z^n$. 

Since $T$ is $\Z$-unimodular, $\varphi_T$ turns out to be {\em bijective}; that is,
$\varphi_T$ is an automorphism of the $\R$-algebra $\ash$.
As a consequence, an {\em ideal} $\a\subseteq\ash$ is mapped to another {\em ideal}
$\varphi_T(\a)$. The map $\varphi_T$ can be extended to a map from $\ash^\ttw$ to itself:
\begin{equation}\label{eq:def:maphi}
\begin{array}{rcl}
\maphi:\ash^\ttw \!&\rightarrow&\! \ash^\ttw\\
\left[\begin{array}{cccc}\!f_1(\uld{}) \!&\! f_2(\uld{}) \!&\! \cdots \!&\! f_\ttw(\uld{})\end{array}\!
\right] \!&\mapsto&\! \left[\begin{array}{cccc}\!\varphi_T(f_1(\uld{})) \!&\! \varphi_T(f_2(\uld{})) \!&\! \cdots \!&\! \varphi_T(f_\ttw(\uld{}))\end{array}\!\right].
\end{array}
\end{equation}
The map $\maphi$ is an $\ash$-module morphism via the automorphism $\varphi_T$,
{\em i.e.}, for $r(\uld{})\in\ash^\ttw$ and $f(\uld{})\in\ash$,
$$
\maphi(f(\uld{})r(\uld{}))=\varphi_T(f(\uld{}))\maphi(r(\uld{})).
$$
The bijective property of $\varphi_T$ extends to the module case: as a result,
$\maphi(\Rmod)$, the image
of a {\em submodule} $\Rmod\subseteq\ash^\ttw$ under $\maphi$, too, is a {\em submodule}.

\subsection{Transforming a single polynomial}
We start off by normalizing a single polynomial. In Lemma \ref{lem:nnl1} below we show that given any
Laurent polynomial $f$ there exists a coordinate transformation matrix $T$ such that
the image of $f$ under the corresponding automorphism $\varphi_T$ has a special
structure: $\varphi_T(f)$ is a linear combination of distinct powers of $\Dx{n}$ with
all coefficients being units in the smaller ring $\ashd{n-1}:=\R[\Dx{1},\Dx{1}^{-1},
\ldots,\Dx{n-1},\Dx{n-1}^{-1}]$. This result is not new;
a proof of this result is given in \cite{park}, where it was used in the context of
designing inverse filters for $n$-D discrete filters. We provide a proof of Lemma \ref{lem:nnl1} here
to make this exposition self-contained.

\newcommand{\dotprod}[1]{\langle\widetilde{\bf t},{#1}\rangle}
\BL\label{lem:nnl1}
Let $f\in\ash$ be given by
$$
f=\sum_{\nu\in\Z^n}\alpha_\nu\uld{}^\nu,~~\mbox{where }\alpha_\nu\in\R,\mbox{ and only finitely many }\alpha_\nu\ne 0.
$$
Then there exists a co-ordinate transformation $T\in\Z^{n\times n}$ such that
under the corresponding homomorphism $\varphi_T:\ash\rightarrow\ash$, we have
\begin{equation}\label{eq:nnl1}
\varphi_T(f)=\left(\sum_{i=0}^{\delta}u_i\Dx{n}^i\right)u,
\end{equation}
where $u_i\in\ashd{n-1}$
and $u\in\R[\Dx{n},\Dx{n}^{-1}]$ are units and $\delta$ is some finite positive
integer.
\EL
\BP
In general, there are many $T$'s which will render $f$ in the form of equation
(\ref{eq:nnl1}); we construct one particular $T$ for which the monomials of
$f$ are mapped to monomials having different $\Dx{n}$ powers. 
To this end, let us define $T\in\Z^{n\times n}$ as
$$
T=\left[\begin{array}{cccc}1 & 0 & \cdots & 0\\
                           0 & 1 & \cdots & 0\\
                           \vdots & \vdots & \ddots & \vdots\\
                           t_1 & t_2 & \cdots & 1\end{array}\right],
$$
where $t_1,t_2,\ldots, t_{n-1}\in\Z$.
Clearly, $T$ is unimodular. The homomorphism $\varphi_T$ corresponding
to this $T$ maps a monomial $\uld{}^\nu$ with $\nu={\rm col}(\nu_1,\nu_2,\ldots,\nu_n)\in
\Z^n$ to the following
$$
\varphi_T(\uld{}^\nu)=\uld{}^{T\nu}=
\Dxt{1}^{\nu_1}\Dxt{2}^{\nu_2}\cdots\Dxt{n-1}^{\nu_{n-1}}\Dxt{n}^{t_1\nu_1+t_2\nu_2+
\cdots+\nu_n}.
$$
Let us define ${\bf t}:={\rm col}(t_1,t_2,\ldots,t_{n-1})\in\Z^{n-1}$, and $\widetilde{\bf t}:={\rm col}({\bf t},1)
\in\Z^n$. Note that the power of $\Dx{n}$ in
$\varphi_T(\uld{}^{\nu})$ is then given by the usual dot product $\dotprod{\nu}$.
Let us also define the following finite subset of $\Z^n$.
$$
{\rm supp}(f):=\{\nu\in\Z^n~|~\alpha_\nu\ne 0\}.
$$ 
We now claim that ${\bf t}$ can be chosen such that the $\Dx{n}$-degrees of $\varphi_T(\uld{}^\nu)$ for
$\nu\in{\rm supp}(f)$ are all different from each other, {\em i.e.}, the following set
\begin{equation}\label{eq:claim1}
\left\{{\bf t}\in\Z^{n-1}~\vline~\dotprod{\nu}\ne\dotprod{\nu'}
\mbox{ for all }\nu,\nu'\in{\rm supp}(f), \nu\ne\nu'\right\}\ne\emptyset.
\end{equation}
To verify this claim let us choose two distinct elements from ${\rm supp}(f)$,
say $\nu={\rm col}(\nu_1,\nu_2,\ldots,\nu_n)$ and $\nu'={\rm col}
(\nu'_1,\nu'_2,\ldots,\nu'_n)$ and consider the following equation
\begin{equation}\label{eq:aux7}
{\bf t}^{\rm T}\left[\begin{array}{c}\nu_1-\nu'_1\\
\nu_2-\nu'_2\\\vdots\\\nu_{n-1}-\nu'_{n-1}\end{array}\right]=\nu'_n-\nu_n.
\end{equation}
We consider solutions to this equation (\ref{eq:aux7}) with ${\bf t}\in\R^{n-1}$,
{\em i.e.},
$$
\mathcal{S}_{\nu,\nu'}=\left\{{\bf t}\in\R^{n-1}~|~\mbox{equation (\ref{eq:aux7}) is
satisfied}\right\}.
$$
The solution set $\mathcal{S}_{\nu,\nu'}$, if nonempty, is an affine hyperplane of
dimension $(n-2)$.
This follows from the following argument. Since
$\nu\ne\nu'$, the vector ${\rm col}((\nu_1-\nu'_1),
(\nu_2-\nu'_2),\ldots,(\nu_{n-1}-\nu'_{n-1}))$ and the number
$(\nu'_n-\nu_n)$ cannot both be equal to zero simultaneously. Hence,
$\mathcal{S}_{\nu,\nu'}$ is either
empty or an affine hyperplane of dimension $(n-2)$. Now consider all pairs $\nu,\nu'\in
{\rm supp}(f)$ and the corresponding $\mathcal{S}_{\nu,\nu'}$'s. Clearly,
there are only finitely many such pairs (in fact, their number is
$\frac{|{\rm supp}(f)|(|{\rm supp}(f)|-1)}{2}$). Consider the following union
$$
\mathcal{S}:=\bigcup_{\nu,\nu'\in{\rm supp}(f),\nu\ne\nu'}\mathcal{S}_{\nu,\nu'}.
$$
Now observe that the set under consideration in equation (\ref{eq:claim1}) is related with
the above defined ${\mathcal{S}}$ as
$$
\left\{{\bf t}\in\Z^{n-1}~|~\dotprod{\nu}\ne\dotprod{\nu'}
\mbox{ for all }\nu,\nu'\in{\rm supp}(f), \nu\ne\nu'\right\}
=(\R^{n-1}\backslash\mathcal{S})\cap\Z^{n-1}.
$$
Since the number of pairs $(\nu,\nu')$ is finite and each of
$\mathcal{S}_{\nu,\nu'}$ is either empty or a hyperplane of dimension $(n-2)$,
the union $\mathcal{S}$ cannot contain all of $\Z^{n-1}$. And therefore,
$(\R^{n-1}\backslash\mathcal{S})\cap\Z^{n-1}\ne\emptyset$, which proves our claim
equation (\ref{eq:claim1}).

It then follows that we can choose a ${\bf t}\in\Z^{n-1}$ such that
$$
\dotprod{\nu}\ne\dotprod{\nu'}
\mbox{ for all }\nu,\nu'\in{\rm supp}(f), \nu\ne\nu'.
$$
In that case, the the nonzero monomials appearing in $\varphi_T(f)$ will have all
distinct powers of $\Dx{n}$. Therefore, $\varphi_T(f)$ can be written as
$$
\varphi_T(f)
=\sum_{\nu\in{\rm supp}(f)}\left(\alpha_\nu\prod_{j=1}^{n-1}
\Dx{j}^{\nu_j}\right)\Dx{n}^{\dotprod{\nu}}.
$$
Note that $\left(\alpha_\nu\prod_{j=1}^{n-1}\Dxt{j}^{\nu_j}\right)$ are all units in $\ashd{n-1}$. Rearranged in
the ascending order of powers of $\Dxt{n}$ we get $u_i$'s
from $\left(\alpha_\nu\prod_{j=1}^{n-1}\Dxt{j}^{\nu_j}\right)$s. Further,
taking out the smallest power of $\Dxt{n}$ and calling it $u$ we get the required expression 
of equation (\ref{eq:nnl1}).
\EP

There is an important consequence of Lemma \ref{lem:nnl1}; we state this as Lemma
\ref{cor:nnl} below. This will be essential for the general case Theorem
\ref{thm:nnl}, where we deal with normalization of a general (non-principal) ideal. 

\BL\label{cor:nnl}
Let $\{0\}\ne\a\subseteq\ash$ be an ideal. Then there exists a coordinate transformation
matrix $T\in\Z^{n\times n}$, inducing $\varphi_T:\ash\rightarrow
\ash$, such that
\begin{itemize}
\item $\ash/\varphi_T(\a)$ contains $\ashd{n-1}/(\varphi_T(\a)\cap\ashd{n-1})$ as
a subring, and 
\item $\ash/\varphi_T(\a)$ is a finitely generated module over $\ashd{n-1}/(\varphi_T(\a)
\cap\ashd{n-1})$.
\end{itemize}
\EL
\BP
First note that for an ideal $\a\subseteq\ash$, the quotient ring $\ash/\a$ can be
written as
$$
\ash/\a=\frac{\ashd{n-1}}{\a\cap\ashd{n-1}}[\overline{\Dx{n}},\overline{\Dx{n}}^{-1}].
$$
This is the algebra generated by $\overline{\Dx{n}}$ and $\left(\overline{\Dx{n}}\right)^{-1}$ over
$\ashd{n-1}/(\a\cap\ashd{n-1})$.
As we have already seen, $\varphi_T$ takes an ideal
$\a\subseteq\ash$ to an ideal $\bid=\varphi_T(\a)\subseteq\ash$. 
Suppose $f\in\a$ is any nonzero Laurent polynomial. 
By Lemma \ref{lem:nnl1}, there exists a coordinate
transformation matrix $T$ such that $\varphi_T(f)$ has the form given by
equation (\ref{eq:nnl1}). Since $u$ in the right-hand-side of equation (\ref{eq:nnl1})
is a unit in $\ash$, it follows that 
$$
\varphi_T(f)=\left(\sum_{i=0}^{\delta}u_i\Dx{n}^i\right)u\in\bid~~~\Leftrightarrow~~~
\sum_{i=0}^{\delta}u_i\Dx{n}^i\in\bid.
$$
Further, since $\{u_i\}$ are all units, the coefficient of the highest power of $\Dx{n}$, {\em i.e.} $u_\delta$, can be made $1$, and the
resulting polynomial will still belong to the same ideal. That is,
$$
u_\delta\left(\Dx{n}^\delta + \sum_{i=0}^{\delta-1}u'_i\Dx{n}^i\right)\in\bid~~~\Leftrightarrow~~~
\Dx{n}^\delta + \sum_{i=0}^{\delta-1}u'_i\Dx{n}^i\in\bid,
$$
where $u'_i=u_iu_\delta^{-1}\in\ash_{n-1}$. In other words, $\overline{\Dx{n}}$
satisfies a monic polynomial equation in $\ash/\bid$ with coefficients coming from
$\ashd{n-1}/(\bid\cap\ashd{n-1})$:
\begin{equation}\label{eq:monic_poly}
\left(\overline{\Dx{n}}\right)^\delta+\sum_{i=0}^{\delta-1}\overline{u'_i}\left(\overline{\Dx{n}}\right)^i=0.
\end{equation}

Now, suppose we have a Laurent polynomial, say $g\in\ash$, which can be written as a polynomial
with non-negative powers in $\Dx{n}$ with coefficients from $\ashd{n-1}$. That is,
$g\in\ashd{n-1}[\Dx{n}]$. Then, we can carry out Euclidean division algorithm by the monic
polynomial appearing in equation (\ref{eq:monic_poly}) above
to obtain 
$$
g=\left(\Dx{n}^\delta + \sum_{i=0}^{\delta-1}u'_i\Dx{n}^i\right)q+\sum_{i=0}^{\delta-1}
\alpha_i\Dx{n}^i,
$$
where $q\in\ashd{n-1}[\Dx{n}]$ and $\{\alpha_0,\alpha_1,\ldots,\alpha_{\delta-1}\}
\subseteq\ashd{n-1}$. Equation (\ref{eq:monic_poly}) then implies
$$
\overline{g}=\sum_{i=0}^{\delta-1}\overline{\alpha_i}\left(\overline{\Dx{n}}\right)^i.
$$
Note that the coefficients above are all in $\ashd{n-1}/(\bid\cap\ashd{n-1})$, and
$\ol{g}\in\frac{\ashd{n-1}}{\bid\cap\ashd{n-1}}[\overline{\Dx{n}}]$. Thus,
every element in $\frac{\ashd{n-1}}{\bid\cap\ashd{n-1}}[\overline{\Dx{n}}]$ can be written
as a finite linear combination of $\{\overline{1},\overline{\Dx{n}},\ldots,
\left(\overline{\Dx{n}}\right)^{\delta-1}\}$ with coefficients from $\ashd{n-1}/(\bid\cap\ashd{n-1})$.
In other words, $\frac{\ashd{n-1}}{\bid\cap\ashd{n-1}}[\overline{\Dx{n}}]$ is
an $\left(\frac{\ashd{n-1}}{\bid\cap\ashd{n-1}}\right)$-module generated by
$\{\overline{1},\overline{\Dx{n}},\ldots,\left(\overline{\Dx{n}}\right)^{\delta-1}\}$.

Further, multiplying both sides of equation (\ref{eq:monic_poly}) by
$\left(\overline{\Dx{n}}\right)^{-1}$ we get
\begin{equation}\label{eq:monic_poly2}
\left(\overline{\Dx{n}}\right)^{\delta-1}+\sum_{i=1}^{\delta-1}\overline{u'_i}
\left(\overline{\Dx{n}}\right)^{i-1}+
\overline{u'_0}\left(\overline{\Dx{n}}\right)^{{-1}}=0.
\end{equation}
Since $u'_0$ was obtained by multiplying two units, namely $u_0$ and $u_\delta^{-1}$, it follows
that $\overline{u'_0}$ is a unit itself. This is because units in $\ash$ go to units in $\ash/\bid$.
Therefore equation (\ref{eq:monic_poly2}) can be rewritten as
\begin{equation}\label{eq:sigma_inv}
\left(\overline{\Dx{n}}\right)^{-1}=-\left(\overline{{u'}_0}\right)^{-1}\left(\overline{\Dx{n}}\right)^{\delta-1}
-\sum_{i=1}^{\delta-1}\left(\overline{{u'}_0}\right)^{-1}\overline{u'_i}\left(\overline{\Dx{n}}\right)^{i-1}.
\end{equation}
By taking positive powers in both sides of equation
(\ref{eq:sigma_inv}), we get that for every $i\in\N$, the monomial $\left(\ol{\Dx{n}}\right)^{-i}$ can be
written as a polynomial in $\frac{\ashd{n-1}}{\bid\cap\ashd{n-1}}[\overline{\Dx{n}}]$. Therefore,
every polynomial in $\frac{\ashd{n-1}}{\bid\cap\ashd{n-1}}[\overline{\Dx{n}}^{-1}]$ too can be
written as a polynomial in $\frac{\ashd{n-1}}{\bid\cap\ashd{n-1}}[\overline{\Dx{n}}]$. In other
words,
$$
\frac{\ashd{n-1}}{\bid\cap\ashd{n-1}}[\overline{\Dx{n}}^{-1}]=
\frac{\ashd{n-1}}{\bid\cap\ashd{n-1}}[\overline{\Dx{n}}].
$$
Hence
\begin{equation}\label{eq:fg}
\ash/\bid=
\frac{\ashd{n-1}}{\bid\cap\ashd{n-1}}[\overline{\Dx{n}},\overline{\Dx{n}}^{-1}]=
\frac{\ashd{n-1}}{\bid\cap\ashd{n-1}}[\overline{\Dx{n}}].
\end{equation}
We have already shown that $\frac{\ashd{n-1}}{\bid\cap\ashd{n-1}}[\overline{\Dx{n}}]$
is a finitely generated $\left(\frac{\ashd{n-1}}{\bid\cap\ashd{n-1}}\right)$-module. 
Therefore, by equation (\ref{eq:fg}) above, we get that $\ash/\bid$ too is a finitely
generated $\left(\frac{\ashd{n-1}}{\bid\cap\ashd{n-1}}\right)$-module. Also note that $\ash/\bid$
contains $\ashd{n-1}/(\bid\cap\ashd{n-1})$ as a subring. This completes the proof.
\EP

\BR\label{rmk:integral}
In the situation of Lemma \ref{cor:nnl}, it follows that every element in $\ash/\varphi_T(\a)$
satisfies a {\em monic polynomial} equation over the subring $\ashd{n-1}/(\varphi_T(\a)
\cap\ashd{n-1})$ (see \cite[Proposition 5.1]{atiyah}). In this case $\ash/\varphi_T(\a)$ is said
to be an {\em integral extension} over the subring $\ashd{n-1}/(\varphi_T(\a)\cap\ashd{n-1})$, or equivalently,
the map $\ashd{n-1}/(\varphi_T(\a)\cap\ashd{n-1})\rightarrow\ash/\varphi_T(\a)$ is {\em injective and
integral}. The converse of this, that is injective and integral implying the bigger ring being
finitely generated module over the smaller ring, is not true in general. However, if the bigger ring
is a finitely generated algebra over the smaller ring, then the converse is also true
(\cite[Corollary 5.2]{atiyah}). In our case, for every $1\leqslant d
\leqslant n$ and an arbitrary ideal $\a\subseteq\ash$, the quotient ring $\ash/\a$ is the algebra
generated by $\{\ol{\Dx{d+1}},\left(\ol{\Dx{d+1}}\right)^{-1},\ldots,\ol{\Dx{n}},\left(\ol{\Dx{n}}\right)^{-1}\}$
over $\ashd{d}/(\a\cap\ashd{d})$. Thus, $\ash/\a$ is a finitely generated algebra over
$\ashd{d}/(\a\cap\ashd{d})$. Therefore, in all the cases
that follow from now, the following two statements are equivalent:
\begin{enumerate}
\item A ring map $\ashd{d}\rightarrow\ash/\a$ is injective and integral.
\item The bigger ring $\ash/\a$ is a finitely generated faithful module over the smaller ring
$\ashd{d}$.
\end{enumerate}
We use the two situations interchangeably in the sequel.  
\ER

\subsection{The general theorem}
Lemma \ref{cor:nnl} above shows that for a nonzero ideal $\a$, there exists a coordinate
transformation matrix $T\in\Z^{n\times n}$,
such that ${\ashd{n-1}}/(\varphi_T(\a)\cap\ashd{n-1})\rightarrow\ash/\varphi_T(\a)$
is {\em injective and integral}. Note that the smaller ring $\ashd{n-1}/(\varphi_T(\a)
\cap\ashd{n-1})$ may not be a Laurent polynomial ring because the intersection ideal
$\varphi_T(\a)\cap\ashd{n-1}$ may not be the zero ideal. However, it is desirable to
have the smaller ring to be a Laurent polynomial ring itself.
Theorem \ref{thm:nnl} shows how this can be done with a suitable coordinate
transformation $T$. We shall require the following Lemma \ref{lem:step:down} to prove
Theorem \ref{thm:nnl}. As done in \cite{cocoa2} for the conventional Noether's
Normalization Lemma, we also follow an induction strategy for the proof of Theorem \ref{thm:nnl}.
While Lemma \ref{cor:nnl} provides the first step of the induction, the inductive step is
what Lemma \ref{lem:step:down} deals with. The lemma tells us that if we already have an injective
and integral ring map $\ashd{d}/(\a\cap\ashd{d})\rightarrow\ash/\a$ such that
$\a\cap\ashd{d}\ne\{0\}$,
then by using a coordinate transformation $T$ we can go down one step further and have
$\ashd{d-1}/(\varphi_T(\a)\cap\ashd{d-1})\rightarrow \ash/\varphi_T(\a)$ injective and integral.

\BL\label{lem:step:down}
Let $\{0\}\ne\a\subseteq\ash$ be an ideal and $1\leqslant d \leqslant n$ an integer
such that $\ashd{d}/(\a\cap\ashd{d})\rightarrow\ash/\a$ is injective and integral. 
Suppose that $\a\cap\ashd{d}\ne\{0\}$.
Then there exists
$T\in\Z^{n\times n}$, inducing $\varphi_T:\ash\ni\uld{}^\nu\mapsto\uld{}^{T\nu}\in
\ash$ for all $\nu\in\Z^n$, such that $\ashd{d-1}/(\varphi_T(\a)\cap\ashd{d-1})\rightarrow
\ash/\varphi_T(\a)$ is injective and integral.
\EL
\BP
We first define $\widetilde{\a}:=\a\cap\ashd{d}$ an ideal in $\ashd{d}$, and then apply
Lemma \ref{cor:nnl} to $\widetilde{\a}$ and $\ashd{d}$. Since $\widetilde{\a}\ne\{0\}$,
by Lemma \ref{cor:nnl},
there exists a unimodular $\widetilde{T}\in\Z^{d\times d}$ inducing the automorphism
$\varphi_{\widetilde{T}}:\ashd{d}\rightarrow\ashd{d}$ such that
\begin{equation}\label{eq:tilde1}
\ashd{d-1}/(\varphi_{\widetilde{T}}(\widetilde{\a})\cap\ashd{d-1})\rightarrow
\ashd{d}/\varphi_{\widetilde{T}}(\widetilde{\a})
\end{equation}
is injective and integral.

With this $\widetilde{T}$ we now define our desired coordinate transformation on $\Z^n$ as
$T:=\left[\begin{array}{cc}
\widetilde{T} & 0\\0 & I_{n-d}\end{array}\right]\in\Z^{n\times n}$. Note that the
restriction of $T$ on the first $d$ coordinate axes is $\widetilde{T}$ and that on the rest
of the coordinate axes is identity. Therefore, the corresponding $\varphi_T:\ash\rightarrow\ash$
maps elements from $\ashd{d}$ to itself. In fact, $\varphi_T|_{\ashd{d}}=\varphi_{\widetilde{T}}$.
Since $\varphi_{\widetilde{T}}$ is an automorphism of $\ashd{d}$, we get that
\begin{equation}\label{eq:thus}
\varphi_T(\ashd{d})
=\varphi_{\widetilde{T}}(\ashd{d})=\ashd{d}.
\end{equation}
Now, since action of $\varphi_T$ restricted to $\ashd{d}$ matches that of $\varphi_{\widetilde{T}}$
on $\ashd{d}$ we get $\varphi_{\widetilde{T}}(\widetilde{\a})=\varphi_T(\widetilde{\a})$.
Note that, since $\varphi_T$ is a bijection, we have $\varphi_T(\a\cap\ashd{d})=
\varphi_T(\a)\cap\varphi_T(\ashd{d})$. So, we can write 
$$
\varphi_T(\widetilde{\a})=\varphi_T(\a)\cap
\varphi_T(\ashd{d})=\varphi_T(\a)\cap\ashd{d}\mbox{ (by equation (\ref{eq:thus}))}.
$$
Since $\ashd{d-1} \subseteq\ashd{d}$, we get from the above equation  
$$
\varphi_{\widetilde{T}}(\widetilde{\a})\cap\ashd{d-1}=\varphi_T(\a)\cap\ashd{d-1}.
$$
Thus equation (\ref{eq:tilde1}) can be rewritten as:
\begin{equation}\label{eq:tilde2}
\ashd{d-1}/(\varphi_{{T}}({\a})\cap\ashd{d-1})\rightarrow
\ashd{d}/(\varphi_{{T}}({\a})\cap\ashd{d})
\end{equation}
is injective and integral.

We now claim that $\ashd{d}/(\varphi_{{T}}({\a})\cap\ashd{d})\rightarrow\ash/\varphi_T(\a)$
is injective and integral. In order to see this, first note that, like before, $\ash/\varphi_T(\a)$
can be written as a finitely generated algebra over the smaller ring $\ashd{d}/(\varphi_{{T}}({\a})
\cap\ashd{d})$ as
$$
\ash/\varphi_T(\a)=\frac{\ashd{d}}{\varphi_{{T}}({\a})\cap\ashd{d}}[\overline{\Dx{d+1}},
\left(\overline{\Dx{d+1}}\right)^{-1},\ldots,\overline{\Dx{n}},\left(\overline{\Dx{n}}\right)^{-1}].
$$
So, in order to prove our claim, it is enough to show that each of the generators of the algebra
written above is integral over the smaller ring. That is, for all $d+1\leqslant i \leqslant n$,
$\overline{\Dx{i}}$ and $\left(\overline{\Dx{i}}\right)^{-1}$ satisfy monic polynomial equations over 
$\ashd{d}/(\varphi_{{T}}({\a})\cap\ashd{d})$. Now, by assumption, $\ash/\a$ is integral over $\ashd{d}/(\a\cap
\ashd{d})$. This means for $d+1\leqslant i \leqslant n$, $\ol{\Dx{i}}$ and $\left(\ol{\Dx{i}}\right)^{-1}$
satisfy monic polynomial equations over $\ashd{d}/(\a\cap\ashd{d})$. That is,
for $d+1\leqslant i \leqslant n$, there exist positive integers
$\delta_i$ and $\gamma_i$, and $a_{i0},a_{i1},\ldots,a_{i\delta_i-1}\in
\ashd{d}$ and $b_{i0},b_{i1},\ldots,b_{i\gamma_i-1}\in\ashd{d}$ such that we have
\begin{equation}\label{eq:int1}
\left.\begin{array}{ccc}
\Dx{i}^{\delta_i}+a_{i\delta_i-1}\Dx{i}^{\delta_i-1}+\cdots+a_{i1}\Dx{i}+a_{i0}&\in&\a\\
(\Dx{i}^{-1})^{\gamma_i}+b_{i\gamma_i-1}(\Dx{i}^{-1})^{\gamma_i-1}+\cdots+b_{i1}\Dx{i}^{-1}+b_{i0}&\in&\a
\end{array}\right\}
\end{equation}
We then apply $\varphi_T$ to equation (\ref{eq:int1}) and use the fact that for
$d+1\leqslant i \leqslant n$ we have $\varphi_T(\Dx{i})=\Dx{i}$ and $\varphi_T(\Dx{i}^{-1})=\Dx{i}^{-1}$.
Thus we get
\begin{equation}\label{eq:int2}
\left.\begin{array}{ccc}
\Dx{i}^{\delta_i}+\varphi_T(a_{i\delta_i-1})\Dx{i}^{\delta_i-1}+\cdots+\varphi_T(a_{i1})\Dx{i}+\varphi_T(a_{i0})&\in&\varphi_T(\a)\\
(\Dx{i}^{-1})^{\gamma_i}+\varphi_T(b_{i\gamma_i-1})(\Dx{i}^{-1})^{\gamma_i-1}+\cdots+\varphi_T(b_{i1})\Dx{i}^{-1}+\varphi_T(b_{i0})&\in&\varphi_T(\a)
\end{array}\right\}
\end{equation}
As we have already seen, from the definition of $\varphi_T$ it follows that $\ashd{d}$ is invariant under $\varphi_T$,
that is, $\varphi_T$ maps elements from $\ashd{d}$ to itself. Therefore, 
$\varphi_T(a_{i0}),\varphi_T(a_{i1}),\ldots,\varphi_T(a_{i\delta_i-1})\in\ashd{d}$ and
$\varphi_T(b_{i0}),\varphi_T(b_{i1}),\ldots,\varphi_T(b_{i\gamma_i-1})\in\ashd{d}$.
Thus equation (\ref{eq:int2}) shows that for all $d+1\leqslant i \leqslant n$,
both $\overline{\Dx{i}},\left(\overline{\Dx{i}}\right)^{-1}\in\ash/\varphi_T(\a)$ satisfy monic polynomial
equations with coefficients coming from $\ashd{d}/(\varphi_{{T}}({\a})\cap\ashd{d})$,
which means $\ash/\varphi_T(\a)$ is integral over $\ashd{d}/(\varphi_T(\a)\cap\ashd{d})$. That
$$\ashd{d}/(\varphi_T(\a)\cap\ashd{d})\rightarrow\ash/\varphi_T(\a)$$ is injective is obvious.
Therefore our claim that $\ashd{d}/(\varphi_{{T}}({\a})\cap\ashd{d})\rightarrow\ash/\varphi_T(\a)$
is injective and integral is true.

By combining equation (\ref{eq:tilde2}) with this we get the following two-stage tower of integral
ring extensions:
$$
\ashd{d-1}/(\varphi_{{T}}({\a})\cap\ashd{d-1})\rightarrow
\ashd{d}/(\varphi_{{T}}({\a})\cap\ashd{d})\rightarrow\ash/\varphi_T(\a).
$$
It is a well-known result in commutative algebra that a finite tower of integral ring extensions
is itself integral (see \cite[Corollary 5.4]{atiyah}). Therefore, we get that  
$$
\ashd{d-1}/(\varphi_{{T}}({\a})\cap\ashd{d-1})\rightarrow
\ash/\varphi_T(\a)
$$
is integral. Injectivity of the above map is clear. This completes the proof.
\EP

We are now in a position to state and prove the main result: Theorem \ref{thm:nnl}.
We call Theorem \ref{thm:nnl} the discrete version of Noether's Normalization Lemma, or DNNL in
short.

\begin{figure}[h!]
\centerline{
\resizebox{16cm}{16cm}{
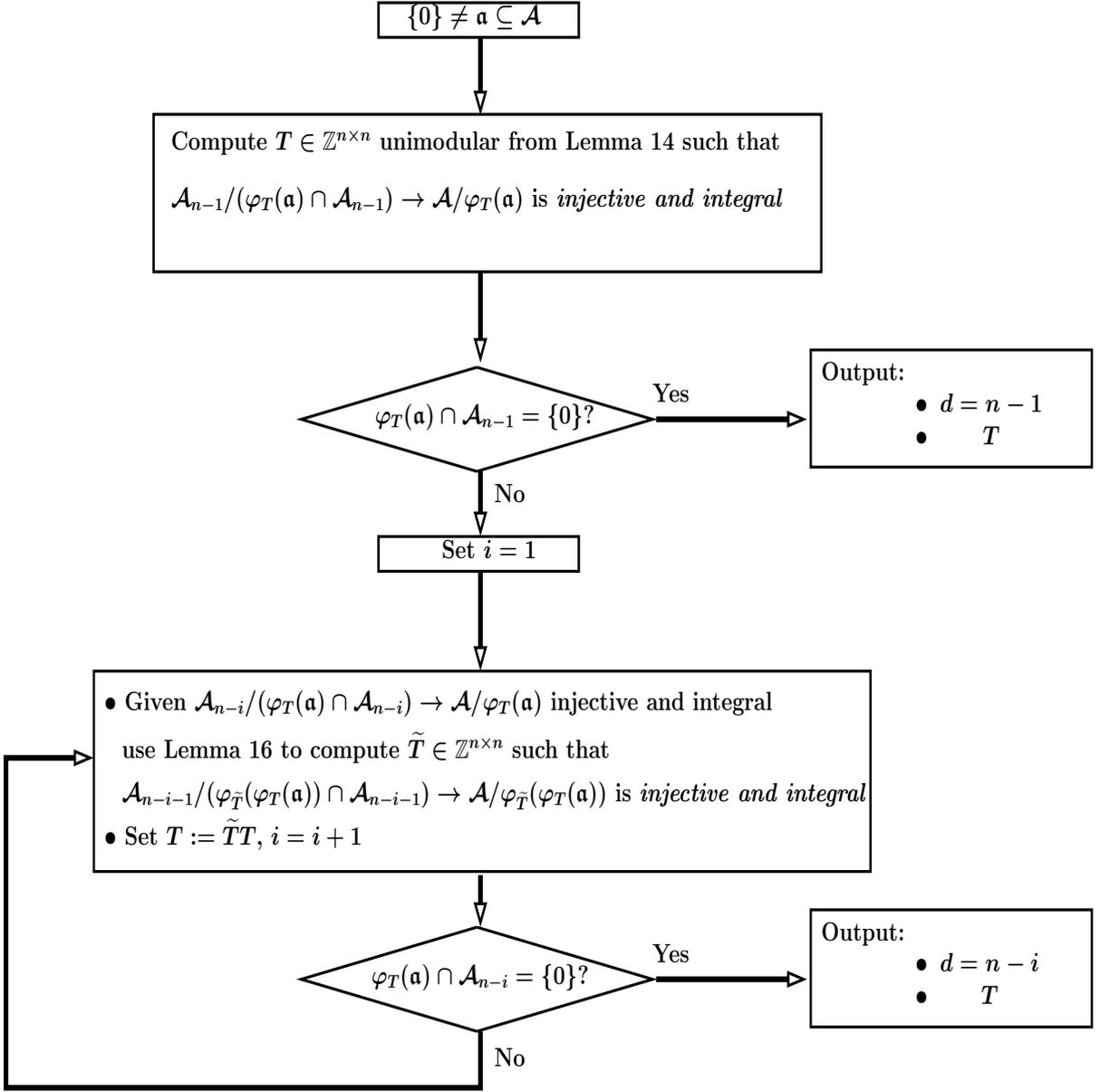
}
}
\caption{Flow-chart used in the proof of Theorem \ref{thm:nnl}}\label{fig:flowchart}
\end{figure}

\BT{\bf (DNNL)}\label{thm:nnl}
Let $\a\subseteq\ash$ be an ideal. Then there exists a coordinate transformation
$T\in\Z^{n \times n}$ inducing the ring map $\varphi_T:\ash\ni\uld{}^\nu\mapsto
\uld{}^{T\nu}\in\ash$ for all $\nu\in\Z^n$, and a positive integer $0\leqslant d \leqslant n$
such that the canonical ring map
$$
\ashd{d}\rightarrow\ash/\varphi_T(\a)
$$
is injective and integral. Equivalently, the following two conditions hold:
\begin{itemize}
\item $\ash/\varphi_T(\a)$ contains $\ashd{d}$ as a subring, and
\item $\ash/\varphi_T(\a)$ is a finitely generated module over $\ashd{d}$.
\end{itemize}
\ET
\BP
As mentioned before, we give a proof by induction on the number $(n-d)$. The first step
of the proof follows from Lemma \ref{cor:nnl}. We may assume that $\a\ne\{0\}$, because
otherwise the theorem trivially holds with $d=n$. Since $\a\ne\{0\}$, by Lemma \ref{cor:nnl}
it follows that there exists $T_1:\Z^{n\times n}$ unimodular such that
$$
\ashd{n-1}/(\varphi_{T_1}(\a)\cap\ashd{n-1})\rightarrow\ash/\varphi_{T_1}(\a)
$$
is injective and integral. The inductive step is given by Lemma \ref{lem:step:down}. If
$\varphi_{T_1}(\a)\cap\ashd{n-1}=\{0\}$ then
the result of the theorem follows with $d=n-1$. On the other hand, if $\varphi_{T_1}(\a)\cap\ashd{n-1}\ne\{0\}$ then by Lemma \ref{lem:step:down} there exists
unimodular $T_2\in\Z^{n\times n}$ such that 
\begin{equation}\label{eq:thm:nnl1}
\ashd{n-2}/(\varphi_{T_2}(\varphi_{T_1}(\a))\cap\ashd{n-2})\rightarrow\ash/\varphi_{T_2}(\varphi_{T_1}(\a))
\end{equation}
is injective and integral.
Note that $\varphi_{T_2}\circ\varphi_{T_1}=\varphi_{T_2T_1}$. Thus equation (\ref{eq:thm:nnl1})
becomes
$$
\ashd{n-2}/(\varphi_{T_2T_1}(\a)\cap\ashd{n-2})\rightarrow\ash/\varphi_{T_2T_1}(\a)
$$
is injective and integral. Next we look at $\varphi_{T_2T_1}(\a)\cap\ashd{n-2}$. If this is the
zero ideal then the theorem holds with $d=n-2$. Otherwise, we once again invoke Lemma
\ref{lem:step:down}, by which we have $T_3\in\Z^{n\times n}$ unimodular, such that
$$
\ashd{n-3}/(\varphi_{T_3T_2T_1}(\a)\cap\ashd{n-3})\rightarrow\ash/\varphi_{T_3T_2T_1}(\a)
$$
is injective and integral. Then we look at $\varphi_{T_3T_2T_1}(\a)\cap\ashd{n-3}$ and depending on
whether this ideal is the zero ideal or not we invoke Lemma \ref{lem:step:down} again or stop.
We illustrate this process of recursively applying Lemma \ref{lem:step:down} in the following
flow-chart Figure \ref{fig:flowchart}.

In the flow-chart, after the first step, if $\varphi_T(\a)\cap\ashd{n-1}\ne\{0\}$
then Lemma \ref{lem:step:down} is applied recursively, and in the process $i$ is incremented at
each step. This recursive process continues as long as $\varphi_T(\a)\cap\ashd{n-i}\ne\{0\}$.
Note that the recursive process is guaranteed to stop after finitely many ($\leqslant n$)
iterations because $i$ is upper-bounded by $n$. When the process stops, we have
the following situation
\begin{equation}\label{eq:thm:nnl2}
\ashd{n-i}/(\varphi_{T}(\a)\cap\ashd{n-i})\rightarrow\ash/\varphi_{T}(\a)
\end{equation}
is injective and integral, plus $\varphi_{T}(\a)\cap\ashd{n-i}=\{0\}$. This means we can rewrite
equation (\ref{eq:thm:nnl2}) as
$$
\ashd{n-i}\rightarrow\ash/\varphi_{T}(\a)
$$
is injective and integral. Defining $d:=n-i$ we get our desired result.
\EP

\BR\label{rmk:krull:dim}
The number $d$ is of special importance to us. This number is equal to the {\em Krull
dimension} of the Laurent polynomial ring $\ashd{d}$. It is a known result in
commutative algebra that the Krull dimensions of two commutative rings $\mathcal{C}_1$
and $\mathcal{C}_2$ are equal if $\mathcal{C}_1\hookrightarrow\mathcal{C}_2$ is
integral (\cite{eisenbud}). It then follows from Theorem \ref{thm:nnl} that
the Krull dimension of $\ash/\varphi_T(\a)$ too is equal to $d$. Yet another fact from
commutative algebra says that Krull dimension remains invariant under isomorphism.
Since $\varphi_T$ is an automorphism of $\ash$, $\ash/\a$ and $\ash/\varphi_T(\a)$ turn
out to be isomorphic to each other as $\R$-algebras. Therefore, the Krull dimension
of $\ash/\a$ is $d$. 
\ER

We illustrate the result of Theorem \ref{thm:nnl} in the following simple example.

\BEx\label{ex:nnl}
Consider the ideal $\a=\langle \Dx{1}\Dx{2}-\Dx{1}-\Dx{2}+1\rangle
\subseteq\R[\Dx{1},\Dx{1}^{-1},\Dx{2},\Dx{2}^{-1}]$. The quotient ring
$\R[\Dx{1},\Dx{1}^{-1},\Dx{2},\Dx{2}^{-1}]/\a$ is not integral over either
$\R[\Dx{1},\Dx{1}^{-1}]$ or $\R[\Dx{2},\Dx{2}^{-1}]$. Theorem \ref{thm:nnl} tells us that with a
coordinate change the quotient ring can be made integral over a smaller Laurent polynomial ring.
Note that ${\rm dim}(\R[\Dx{1},\Dx{1}^{-1},\Dx{2},\Dx{2}^{-1}]/\a)=1$. So the first reduction
in the flow-chart of Figure \ref{fig:flowchart} yields the desired normalization.
By Lemma \ref{lem:nnl1} we can take $T\in\Z^{2\times 2}$ to be
$$
T=\left[\begin{array}{rr}1 & 0 \\2 & 1\end{array}\right].
$$
The corresponding $\varphi_T$ turns out to be $\varphi_T:\Dx{1}\mapsto\Dx{1}\Dx{2}^{2},
\Dx{2}\mapsto\Dx{2}$. This $\varphi_T$ maps the generator of $\a$ to
$$\varphi_T(\Dx{1}\Dx{2}-\Dx{1}-\Dx{2}+1)=\Dx{1}\Dx{2}^3-\Dx{1}\Dx{2}^{2}-\Dx{2}+1
=\Dx{1}\left(\Dx{2}^3-\Dx{2}^2-\Dx{1}^{-1}\Dx{2}+\Dx{1}^{-1}\right).$$
It then follows that
$\varphi_T(\a)=\langle \Dx{2}^3-\Dx{2}^2-\Dx{1}^{-1}\Dx{2}+\Dx{1}^{-1} \rangle$.
Note that the generator of $\varphi_T(\a)$ is a monic polynomial in $\Dx{2}$ with coefficients
coming from $\R[\Dx{1},\Dx{1}^{-1}]$. In other words, $\ol{\Dx{2}}$ in 
$\R[\Dx{1},\Dx{1}^{-1},\Dx{2},\Dx{2}^{-1}]/\varphi_T(\a)$ satisfies a monic polynomial equation
over $\R[\Dx{1},\Dx{1}^{-1}]/(\varphi_T(\a)\cap\R[\Dx{1},\Dx{1}^{-1}])$. Moreover,
$\varphi_T(\a)\cap\R[\Dx{1},\Dx{1}^{-1}]$ is clearly the zero ideal. Thus we get
$\R[\Dx{1},\Dx{1}^{-1}]\rightarrow\R[\Dx{1},\Dx{1}^{-1},\Dx{2},\Dx{2}^{-1}]/\varphi_T(\a)$
is injective and integral.
\EEx

\subsection{DNNL extended to modules}
In this subsection we show how the normalization process can be extended to modules. This extension is essential in order
to apply the normalization process to non-scalar systems.
The extension of DNNL to modules is obtained by two crucial observations, which we elaborate now. In order to proceed further,
we require the following algebraic object: the annihilator ideal of a quotient module $\M$, which is defined as follows:
\begin{equation}\label{eq:def:ann}
{\rm ann}(\M)=\{f(\uld{})\in\ash~|~f(\uld{})m=0\in\M\mbox{ for all }m\in\M\}.
\end{equation}
Note that $\M$ is {\em naturally}\footnote{`Naturally' here means the scalar multiplication for the base rings $\ash$ and
$\ash/{\rm ann}(\M)$ are identical.} a module over $\ash/{\rm ann}(\M)$. Further,
since the free module $\ash^\ttw$ is of finite rank ($=\ttw$) over $\ash$, it is finitely generated as a module
over $\ash$. Now, $\M$, being a quotient of $\ash^\ttw$, is clearly finitely generated over $\ash$. 
Hence, $\M$ is finitely generated over $\ash/{\rm ann}(\M)$ too. This is the first
crucial observation.

The second observation involves the maps induced by a coordinate transformation of $\Z^n$. We write this observation
as a proposition below.

\BPP\label{prop:ann:id}
Let $T\in\Z^{n\times n}$ be a $\Z$-unimodular matrix. Further, let $\varphi_T$ and $\maphi$ be the two maps
induced by $T$ as in equations \eqref{eq:def:phi} and \eqref{eq:def:maphi}, respectively. Suppose $\Rmod\subseteq
\ash^{\ttw}$ is a submodule of the free module $\ash^\ttw$. Then the following holds:
\begin{equation}\label{def:prop:ann:1}
\varphi_T({\rm ann}(\ash^\ttw/\Rmod))={\rm ann}(\ash^\ttw/\maphi(\Rmod)).
\end{equation}    
\EPP
\BP
Note that $f(\uld{})\in{\rm ann}(\ash^\ttw/\maphi(\Rmod))$ if and only if for all $r(\uld{})\in\ash^\ttw$ we have
$f(\uld{})r(\uld{})\in\maphi(\Rmod)$. However, since $\maphi$ is a bijection, for every $r(\uld{})\in\ash^\ttw$
there exists a unique $r_1(\uld{})\in\ash^\ttw$ such that $\maphi(r_1(\uld{}))=r(\uld{})$. Therefore, we can write
$f(\uld{})\in{\rm ann}(\ash^\ttw/\maphi(\Rmod))$ if and
only if for all $r_1(\uld{})\in\ash^\ttw$ we have $f(\uld{})\maphi(r_1(\uld{}))\in\maphi(\Rmod)$. Recall that
$\varphi_T$ is an automorphism of $\ash$. Therefore, there is a unique $g(\uld{}):=\varphi_T^{-1}(f(\uld{}))$.
Then we can write
$$
f(\uld{})\maphi(r_1(\uld{}))=\varphi_T(g(\uld{}))\maphi(r_1(\uld{}))
=\maphi(g(\uld{})r_1(\uld{}))\in\maphi(\Rmod) ~\Leftrightarrow~g(\uld{})\in \ann.
$$
But $g(\uld{})\in{\rm ann}(\ash^\ttw/\Rmod)$ if and only if $f(\uld{})\in\varphi_T({\rm ann}(\ash^\ttw/\Rmod))$.
\EP

We now present the non-scalar or module theoretic version of DNNL.

\BT{\bf (Module theoretic DNNL)}\label{thm:mod:dnnl}
Let $\Rmod\subseteq\ash^\ttw$ be a submodule. Suppose the Krull dimension of ${\rm ann}(\ash^\ttw/\Rmod)$ is $d$.
Then there exists a coordinate transformation
$T\in\Z^{n \times n}$ inducing the $\R$-algebra homomorphism $\varphi_T:\ash\ni\uld{}^\nu\mapsto
\uld{}^{T\nu}\in\ash$ for all $\nu\in\Z^n$, and the corresponding $\ash$-module homomorphism $\maphi:\ash^\ttw\rightarrow\ash^\ttw$,
such that the quotient module $\ash^\ttw/\maphi(\Rmod)$
is a finitely generated module over $\ashd{d}$.
\ET
\BP
Let $T\in\Z^{n\times n}$ be an arbitrary $\Z$-unimodular integer matrix.
We have already seen that for any submodule $\Rmod$ of $\ash^\ttw$, the quotient module $\ash^\ttw/\Rmod$
is finitely generated as a module over $\ash/{\rm ann}(\ash^\ttw/\Rmod)$. In particular, since $\maphi(\Rmod)$ too is
a submodule of $\ash^\ttw$ we get that $\ash^\ttw/\maphi(\Rmod)$ is a finitely generated module over
$\ash/{\rm ann}(\ash^\ttw/\maphi(\Rmod))$.
Now, let us define $\a:={\rm ann}(\ash^\ttw/\Rmod)$.
By Proposition \ref{prop:ann:id},
$\varphi_T(\a)={\rm ann}(\ash^\ttw/\maphi(\Rmod))$.
Therefore, $\ash^\ttw/\maphi(\Rmod)$ is in fact a finitely generated module over $\ash/\varphi_T(\a)$.

Since Krull dimension of $\ash/\a$ is equal to $d$, by applying
Theorem \ref{thm:nnl} to $\a$ we get that there exists a $\Z$-unimodular $T\in\Z^{n\times n}$ such that
$\ash/\varphi_T(\a)$ is a finitely generated module over $\ashd{d}$.  
Thus, by combining this with the fact that $\ash^\ttw/\maphi(\Rmod)$
is a finitely generated module over $\ash/\varphi_T(\a)$ we get the following chain of modules
$$
\begin{array}{c}
\ash^\ttw/\maphi(\Rmod)\\
\uparrow\\
\ash/\varphi_T(\a)\\
\uparrow\\
\ashd{d}
\end{array}
$$ 
where each step is finitely generated. Therefore, $\ash^\ttw/\maphi(\Rmod)$ must be finitely generated over $\ashd{d}$.
Indeed, suppose $\{g_1,\ldots,g_{k}\}\subseteq\ash^\ttw/\maphi(\Rmod)$ is a generating set for $\ash^\ttw/\maphi(\Rmod)$
as a module over $\ash/\varphi_T(\a)$, and $\{h_1,\ldots,h_{\ell}\}\subseteq
\ash/\varphi_T(\a)$ is a generating set for $\ash/\varphi_T(\a)$ as a module
over $\ashd{d}$. Then
$$
\{g_i h_j~|~1\leqslant i \leqslant k,~1\leqslant j \leqslant \ell\}
$$ 
is a generating set for $\ash^\ttw/\maphi(\Rmod)$ over $\ashd{d}$.
\EP

We illustrate this result in the following example.
\BEx
Consider the submodule $\Rmod$ of $\R[\Dx{1},\Dx{1}^{-1},\Dx{2},\Dx{2}^{-1}]^2$ generated by the rows of the matrix
$$
R(\Dx{})=\left[\begin{matrix}\Dx{1}-1 & 2\\1 & \Dx{2}-1\end{matrix}\right].
$$
It can be checked that the the quotient module $\M=\R[\Dx{1},\Dx{1}^{-1},\Dx{2},\Dx{2}^{-1}]^2/\Rmod$ is not finitely generated
over $\R[\Dx{1},\Dx{1}^{-1}]$ or $\R[\Dx{1},\Dx{1}^{-1}]$, neither it is a finite dimensional vector space over $\R$.
The annihilator ideal of the quotient module in this case
turns out to be $\a=\langle\Dx{1}\Dx{2}-\Dx{2}-\Dx{1}-1\rangle$. Take the unimodular matrix $T=\left[\begin{smallmatrix}1 & 0\\
2 & 1\end{smallmatrix}\right]$ for coordinate transformation. Then we get that under the ring map $\varphi_T$ the ideal $\a$
is mapped to the ideal $\varphi_T(\a)=\langle\Dx{1}\Dx{2}^3-\Dx{2}-\Dx{1}\Dx{2}^2-1\rangle$. Note that
$\R[\Dx{1},\Dx{1}^{-1},\Dx{2},\Dx{2}^{-1}]/\varphi_T(\a)$ is a finitely generated module over $\R[\Dx{1},\Dx{1}^{-1}]$.

The extension of $\varphi_T$ to $\R[\Dx{1},\Dx{1}^{-1},\Dx{2},\Dx{2}^{-1}]^2$, that is, $\maphi$ maps $\Rmod$ to the row-span
of the following matrix:
$$
\widetilde{R}(\Dx{})=\left[\begin{matrix}\Dx{1}\Dx{2}^2 -1 & 2\\ 1 & \Dx{2}-1\end{matrix}\right].
$$
Note that the annihilator of $\R[\Dx{1},\Dx{1}^{-1},\Dx{2},\Dx{2}^{-1}]^2/\maphi(\Rmod)$ is given by $\varphi_T(\a)$. It then
follows that $\R[\Dx{1},\Dx{1}^{-1},\Dx{2},\Dx{2}^{-1}]^2/\maphi(\Rmod)$ is finitely generated as a module over
$\R[\Dx{1},\Dx{1}^{-1}]$. Indeed, it can be checked that,
$$
\left\{\left[\begin{matrix}1\\0\end{matrix}\right]^{\rm T},\left[\begin{matrix}0\\1\end{matrix}\right]^{\rm T},
\left[\begin{matrix}\Dx{1}\\0\end{matrix}\right]^{\rm T},\left[\begin{matrix}0\\\Dx{1}\end{matrix}\right]^{\rm T},
\left[\begin{matrix}\Dx{1}^2\\0\end{matrix}\right]^{\rm T},\left[\begin{matrix}0\\\Dx{1}^2\end{matrix}\right]^{\rm T}
\right\}
$$ 
is a generating set for $\R[\Dx{1},\Dx{1}^{-1},\Dx{2},\Dx{2}^{-1}]^2/\maphi(\Rmod)$ as a module over $\R[\Dx{1},\Dx{1}^{-1}]$.
\EEx

\section{Representation formula for general autonomous systems}\label{sec:rep:gen}
In this section we make use of Theorem \ref{thm:mod:dnnl} to obtain representation formula for general
autonomous systems. This clearly is indicative of the fact that this representation formula that we are
going to  propose will involve a coordinate transformation on the indexing set $\Z^n$. Therefore,
we first describe the effects of a coordinate transformation of $\Z^n$ on the behavior of an $n$-D system.

Recall that a coordinate trasnformation of $\Z^n$ is given by a $\Z$-unimodular matrix $T\in\Z^{n\times n}$
as
$$
\begin{array}{rcc}
T:\Z^n & \rightarrow & \Z^n\\
{\rm col}(\nu_1,\nu_2,\ldots,\nu_n)=:\nu & \mapsto & T\nu.
\end{array}
$$ 
Now note that this $T$ induces the following $\R$-linear map from $\disww{n}{\ttw}$ to itself:

\begin{equation}\label{eq:def:pullback}
\begin{array}{rcc}
\Phi_T:\disww{n}{\ttw} &\rightarrow& \disww{n}{\ttw}\\
w(\nu) &\mapsto& w(T\nu),
\end{array}
\end{equation}
This map $\Phi_T$ is often called the {\em pull-back} of the coordinate transformation $T$. The pull-back
$\Phi_T$ is a bijection because $T$ is $\Z$-unimodular. In fact, it is easy to see that
$$
\Phi_T^{-1}=\Phi_{T^{-1}}.
$$

Recall the $\R$-algebra homomorphism $\varphi_T:\ash\rightarrow\ash$ induced by $T$, and its extension to
modules $\maphi:\ash^\ttw\times\ash^\ttw$. The following lemma shows how $\maphi$ and the pull-back $\Phi_T$
are interrelated.

\BL\label{lem:change_of_coord}
Let $v,w\in\disww{n}{\ttw}$ be related by $v=\Phi_T(w)$. Then for $r(\uld{})\in\ash^\ttw$ we have
$$
r(\uld{})v=\Phi_T(\maphi(r(\uld{}))w).
$$
\EL

\BP
It is enough to prove for the scalar case: for $v,w\in\R^{\Z^n}$ and $f(\uld{})\in\ash$
we must have
\begin{equation}\label{eq:cod}
f(\uld{})v=\Phi_T(\varphi_T(f(\uld{}))w).
\end{equation}
This is because $\maphi$ is defined as $\varphi_T$ applied element-wise, and $\Phi_T$ is linear.
Further, note that every Laurent polynomial is a finite $\R$-linear combination of monomials.
Therefore, once again by linearity of $\Phi_T$, it is enough to prove equation (\ref{eq:cod})
for $\ash\ni f(\uld{})=\uld{}^\nu$, where $\nu\in\Z^n$ is arbitrary.

Now, $v,w\in\R^{\Z^n}$ are assumed to be related by $v=\Phi_T(w)$. First note that for
$\nu,\nu'\in\Z^n$,
the action of the monomial $\uld{}^{\nu}$ on $v$ is given by the following
\begin{eqnarray*}
(\uld{}^{\nu}v)(\nu')=v(\nu'+\nu)=\Phi_T(w)(\nu'+\nu)=w(T(\nu'+\nu))&=&(\uld{}^{T\nu}w)(T\nu')\\
&=&\Phi_T(\varphi_T(\uld{}^{\nu})w)(\nu').
\end{eqnarray*}
Since $\nu'$ was arbitrary we get for all $\nu\in\Z^n$
\begin{equation}\label{eq:cod2}
\uld{}^\nu v=\Phi_T(\varphi_T(\uld{}^\nu)w).
\end{equation}
This is what we wanted to show. 
\EP

With this we now bring out the effect of $T$ on the behavior $\B$ of an $n$-D system. We show that the image
of $\B$ under the pull-back $\Phi_T$ is also a behavior. Not only that, in fact, the equation modules of
$\B$ and $\Phi_T(\B)$ are related by $\maphi$ in a way as shown in Proposition \ref{prop:coord:beh} below.

\BPP\label{prop:coord:beh}
Let $\Rmod\subseteq\ash^\ttw$ be a submodule with behavior $\B(\Rmod)$, and let
$T\in\Z^{n\times n}$ be $\Z$-unimodular. Then we have
\begin{equation}\label{eq:thm:change_of_coord}
\Phi_T\left(\B(\Rmod)\right)=\B(\imaphi(\Rmod)).
\end{equation}
\EPP

\BP
Let $w\in\disww{n}{\ttw}$ be an arbitrary $n$-D trajectory and let $r(\uld{})\in\ash^\ttw$ be an arbitrary row-vector
of polynomial shift operators. Then, by Lemma \ref{lem:change_of_coord}, we have the following identity
\begin{equation}\label{eq:prop:coord:1}
r(\uld{})\left(\Phi_T(w)\right) = \Phi_T\left(\maphi(r(\uld{}))w\right).
\end{equation}
Since $\maphi$ and $\Phi_T$ both are bijective, it follows from equation \eqref{eq:prop:coord:1} that
$w\in\B(\maphi(\Rmod))$ if and only if $\Phi_T(w)\in\B(\Rmod)$. In other words,
$$
\Phi_{T^{-1}}\left(\B(\Rmod)\right)=\B(\maphi(\Rmod)).
$$
Since $T$ is invertible we can start the argument with $T$ substituted by $T^{-1}$, and thus we obtain the desired
equation \eqref{eq:thm:change_of_coord}. 
\EP

Proposition \ref{prop:coord:beh} allows us to give the following definition.

\BD{\bf ($T$-isomorphism of behaviors)}\label{def:t-iso}
Suppose $\B_1,\B_2\in\Lwe$. Then $\B_1$ is said to be $T$-isomorphic to $\B_2$ if there exists a $\Z$-unimodular
matrix $T\in\Z^{n\times n}$ such that
$$
\B_1=\Phi_T(\B_2).
$$
\ED

Proposition \ref{prop:coord:beh} then translates in the language of Definition \ref{def:t-iso} to the following:

\BC\label{cor:coord:beh}
Let $\B_1,\B_2\in\Lwe$ have their equation modules $\Rmod_1,\Rmod_2$, respectively. Then $\B_1$ is $T$-isomorphic
to $\B_2$ if and only if $\Rmod_1=\imaphi(\Rmod_2)$, which in turn is true if and only if $\Rmod_2=\maphi(\Rmod_1)$.
\EC

Combining Corollary \ref{cor:coord:beh} (or equivalently, Proposition \ref{prop:coord:beh}) with the module
theoretic DNNL (Theorem \ref{thm:mod:dnnl}) we get the following important result.

\BT\label{thm:main:t-iso}
Let $\B\in\Lwe$ be an autonomous behavior with equation module $\Rmod\subseteq\ash^\ttw$. Further, let the Krull
dimension of ${\rm ann}\left(\ash^\ttw/\Rmod\right)$ be equal to $d$. Then $\B$ is $T$-isomorphic to a strongly relevant autonomous system
of order $d$.
\ET

\BP
First note that, by Theorem \ref{thm:mod:dnnl}, there exists a $\Z$-unimodular $T\in\Z^{n\times n}$ such that
under the induced module map $\maphi:\ash^\ttw\rightarrow\ash^\ttw$ we have $\ash^\ttw/\maphi(\Rmod)$ a
finitely generated module over $\ashd{d}$. Therefore, $\B(\maphi(\Rmod))$ is a strongly relevant autonomous
behavior of order $d$, by Definition \ref{def:str:rel}. By Corollary \ref{cor:coord:beh}, $\B$ is $T$-isomorphic
to $\B(\maphi(\Rmod))$. Thus $\B$ is $T$-isomorphic to a strongly relevant autonomous system of order $d$.
\EP

Theorem \ref{thm:main:t-iso} enables us to obtain the following representations for a general autonomous system.
These representations follow from the representations of strongly relevant systems presented in Theorem
\ref{thm:main:spl} and Corollary \ref{cor:main:spl}. Recall the first order system of equations given in
Table 1. The representations given below involve these equations. The next result is a consequence of Theorem
\ref{thm:main:spl}.

\BC\label{cor:main:1}
Suppose $\B\in\Lwe$ is an autonomous system with equation module $\Rmod$. Let $d$ be the Krull dimension
of ${\rm ann}\left(\ash^\ttw/\Rmod\right)$. Then, there exist a positive integer $\gamma$ and the following matrices:
\begin{itemize}
\item $X(\uld{1})\in\ashd{d}^{\bullet\times\gamma}$,
\item $A_1(\uld{1}),A_2(\uld{1}),\ldots,A_{n-d}(\uld{1})\in\ashd{d}^{\gamma\times\gamma}$,
\item $C(\uld{1})\in\ashd{d}^{\ttw\times\gamma}$,
\end{itemize}
such that $\B$ admits the following representation:
\begin{equation}
\B=\left\{v\in\disww{n}{\ttw}~ \vline\begin{array}{l}\exists ~\xvec\in\disww{n}{\gamma},w\in\disww{n}{\ttw}
 \mbox{such that equations in Table 1}\\ \mbox{are satisfied, and } v=\Phi_T(w)\end{array}\right\}
\end{equation}
\EC 

The next result is a consequence of Corollary \ref{cor:main:spl}
\BC\label{cor:main:2}
Let $\B\in\Lwe$ be an autonomous system whose equation module is $\Rmod\subseteq\ash^\ttw$. Suppose the Krull
dimension of ${\rm ann}(\ash^\ttw/\Rmod)$ is $d$. Then there exist
\begin{itemize}
\item a positive integer $\gamma$,
\item a matrix $X(\uld{1})\in\ashd{d}^{\bullet\times\gamma}$,
\item a set of $n-d$ square matrices $A_1(\uld{1}),\ldots,A_{n-d}(\uld{1})\in\ashd{d}^{\gamma\times\gamma}$,
\item a matrix $C(\uld{1})\in\ashd{d}^{\ttw\times\gamma}$
\end{itemize}
such that $\B$ admits the following latent variable representation:
\begin{equation}\label{eq:cor:lat}
\B=\left\{v\in\disww{n}{\ttw}~\vline~\begin{array}{l}\exists ~\xvec\in\disww{n}{\gamma}\mbox{ and }w\in\disww{n}{\ttw}
\mbox{ satisfying }\\
\left[\begin{array}{c:c}X(\uld{1}) & 0 \\\hdashline
\Dx{d+1}I_{\gamma}-A_1(\uld{1}) & 0\\
\Dx{d+2}I_{\gamma}-A_2(\uld{1}) & 0\\
\vdots & \vdots \\
\Dx{n}I_{\gamma}-A_{n-d}(\uld{1}) & 0\\\hdashline
-C(\uld{1}) & I_{\ttw}
\end{array}\right]\left[\begin{matrix}\xvec \\ w\end{matrix}\right]=0,~ v=\Phi_T(w)\end{array}
 \right\}
\end{equation}
\EC

With these representations for general autonomous systems, we are now in a position to provide an algorithm for
finiding solutions for general autonomous systems. The algorithm is exactly the same as Algorithm \ref{spl-algo}
with one extra stage at the beginning, where we {\em normalize} the given system to a $T$-isomorphic system that is
strongly relevant.  

\BA{\bf (Solution formula for general autonomous systems)}\label{main-algo}

\begin{tabular}{|l|}
\hline
{\bf (Noether) Normalization}~\vline\\
\hline
{\bf Input:} An autonomous system $\B={\rm ker}(R(\uld{}))$, where $R(\uld{})\in\ash^{\bullet\times\ttw}$.\\
\hline
{\bf Computation:}\\
Define $\Rmod={\rm rowspan}(R(\uld{}))\subseteq\ash^\ttw$.\\
Obtain $\a={\rm ann}(\ash^\ttw/\Rmod)$.\\
Using Theorem \ref{thm:nnl} find out a $\Z$-unimodular $T\in\Z^{n\times n}$ such that under\\
the induced map $\varphi_T:\ash\rightarrow\ash$ we have $\ashd{d}\rightarrow\ash/\varphi_T(\a)$
injective and\\ integral.\\
\hline
{\bf Output:}\\ 
$T$.\\
$\widetilde{\Rmod}:=\maphi(\Rmod)$.\\
\hline
{\bf Regularization}~\vline\\
\hline
{\bf Input:} $\widetilde{\Rmod}\subseteq\ash^\ttw$.\\
\hline
{\bf Computation:}\\
Obtain a generating set $\mathcal{G}:=\{g_1,\ldots,g_{\gamma}\}\subseteq\ash^\ttw/\widetilde{\Rmod}$.\\
Obtain the matrix of relations $X(\uld{1})\in\ashd{d}^{k\times\gamma}$.\\
Obtain the companion matrices $A_1(\uld{1}),\ldots,A_{n-d}(\uld{1})\in\ashd{d}^{\gamma\times\gamma}$ (ensure\\
that these matrices are invertible in $\ashd{d}$).\\
Obtain the output matrix $C(\uld{1})\in\ashd{d}^{\ttw\times\gamma}$.\\
\hline
{\bf Output:} The matrices $X(\uld{1})$, $A_1(\uld{1}),\ldots,A_{n-d}(\uld{1})$, and $C(\uld{1})$.\\
\hline
{\bf Recursion}~\vline\\
\hline
{\bf Input:} {\em The set of initial data}: Solve $X(\uld{1})x=0$ for $x\in\disww{d}{\gamma}$.\\
\hline
{\bf Computation:}\\
{\bf for} $\nu=(\nu_1,\nu_2,\ldots,\nu_n)\in\Z^n$ set\\
~~~~~~~~~$
\widetilde{w}(\nu)=\left(C(\uld{1})\prod_{i=1}^{n-d}A_i(\uld{1})^{\nu_{d+i}}x\right)(\nu_1,\nu_2,\ldots,\nu_d).
$\\
{\bf end}\\
\hline
{\bf Output:} The sequence $\widetilde{w}:=\{\widetilde{w}(\nu)\}_{\nu\in\Z^n}$.\\
\hline
{\bf Renormalization}~\vline\\
\hline
{\bf Input:} The sequence $\widetilde{w}:=\{\widetilde{w}(\nu)\}_{\nu\in\Z^n}$.\\
\hline
{\bf Computation:}\\
{\bf for} $\nu=(\nu_1,\nu_2,\ldots,\nu_n)\in\Z^n$ set\\
~~~~~~~~~$
w(\nu)=\widetilde{w}(T\nu).
$\\
{\bf end}\\
\hline
{\bf Output:} The sequence ${w}:=\{{w}(\nu)\}_{\nu\in\Z^n}$.\\
\hline
\end{tabular}
\EA

\section{The set of allowable initial conditions}\label{sec:allow}
In Algorithms \ref{spl-algo} and \ref{main-algo}, the initial conditions $x\in\disww{d}{\gamma}$ are required to
satisfy the equation $X(\uld{1})x=0$. This is in contrast with the $1$-D situation, where 
one can {\em always} choose $x$ in such a way that the set of initial conditions is free
(see \cite{will:83}). In this section we analyze the set of allowable initial conditions
and address certain important issues, namely, {\em freeness} of the initial conditions, and {\em non-autonomy} of the
state space $\X$.

\subsection{Freeness of the set of allowable initial conditions}
We have already seen that the representation of an autonomous syetem given by Theorem \ref{thm:main:spl}
can be viewed as an $(n-d)$-D flow equation arising out of an $(n-d)$-D first order dynamics on an
infinite dimensional space. We start by identifying this infinite dimensional space.
Suppose $\Rmod\subseteq\ash^\ttw$ is a submodule such that $\ash^\ttw/\Rmod$ is a finitely generated module
over $\ashd{d}$. We fix a generating set $\{{g_1},\ldots,{g_\gamma}\}$
for $\ash^\ttw/\Rmod$ as an $\ashd{d}$-module. Recall equation (\ref{eq:finitely_gen}),
where we identified $\ash^\ttw/\Rmod$ with a quotient module of $\ashd{d}^\gamma$ via the map
$\psi:\ashd{d}^\gamma\rightarrow\ash^\ttw/\Rmod$. Let $X(\uld{1})\in\ashd{d}^{\delta\times\gamma}$ be a
matrix of relations of $\{{g_1},\ldots,{g_\gamma}\}$.
We have defined the following set of vector valued $d$-D trajectories:
\begin{equation}\label{eq:def:X}
\X:=\left\{x\in\disww{d}{\gamma}~|~X(\uld{1})x=0\right\},
\end{equation}
as the state space.

It is clear that $\X$ contains all the allowable initial conditions. Note that $\X$ is in fact a vector valued
$d$-D behavior having $X(\uld{1})$ as a kernel representation
matrix. The submodule corresponding to this behavior
is the row-span of the matrix $X(\uld{1})$ over $\ashd{d}$, which equals ${\rm ker}(\psi)$.
Clearly, this submodule ${\rm ker}(\psi)$ is unique for given set of generators of $\ash^\ttw/\Rmod$ as
a module over $\ashd{d}$.
Thus, by the one-to-one correspondence between behaviors and submodules,
$\X$ gets fixed once a generating set for $\ash^\ttw/\Rmod$ (as an $\ashd{d}$-module) is fixed.
However, even for a fixed $\Rmod$, different generating sets for $\ash^\ttw/\Rmod$ as a $\ashd{d}$-module may exist.
These different generating sets may even have different modules of relations, leading to distinct
$\X$s. But, it can be shown that these different $\X$s are {\em isomorphic} to each other
in the sense of \cite{pillai-wood-rogers}, because they have isomorphic quotient modules.
One important question related to $\X$ is: when
does $\ash^\ttw/\Rmod$ admit a generating set of cardinality $\gamma$ such that the corresponding $\X$
equals the whole
of $\disww{d}{\gamma}$? In the sequel, we refer to this situation by saying $\ash^\ttw/\Rmod$ {\em
admits a free $\X$ of rank $\gamma$}.

It is easy to see that the answer to this question is: if and only if $\ash^\ttw/\Rmod$ is a {\em free
module of rank $\gamma$ over $\ashd{d}$}. This is because then, and only then, $\ash^\ttw/\Rmod$ will admit a
free generating set, i.e., the matrix of relations $X(\uld{1})$ corresponding to the generating set will
be the zero matrix. And consequently,
$\X$ will be the whole of $\disww{d}{\gamma}$. This fact together with Swan's extension
\cite{swan} of Serre conjecture to Laurent polynomial rings give us the following check (Theorem
\ref{thm:check:free}) for determining whether $\ash^\ttw/\Rmod$ admits a free $\X$ or not. 

\BR\label{rmk:zero:lp} 
In Theorem
\ref{thm:check:free} we require the notion of a {\em zero left-prime} Laurent polynomial matrix.
A Laurent polynomial matrix $X(\uld{1})\in\ashd{d}^{k \times \gamma}$ having {full row rank} is said to be
{\em zero left-prime} if the ideal generated by the $(k\times k)$ minors of $X(\uld{1})$ is the
whole
ring $\ashd{d}$. See \cite{zerz,bose} for various different notions of prime-ness and their significances.
One particular property of zero left-prime matrices that will be important for proving Theorem \ref{thm:check:free}
is the following:

\noindent
{\em The row-span of a zero left-prime matrix is a direct-summand of the corresponding free module. And, conversely,
if a submodule $\Rmod$ of the free module $\ash^\ttw$ is a direct-summand then $\Rmod$ admits representation as
the row-span of a matrix that is zero left-prime \cite{zerz,bose}. }
\ER

\BT\label{thm:check:free}
Let $\Rmod\subseteq\ash^\ttw$ be a submodule such that $\ash^\ttw/\Rmod$ is a finitely generated module over
$\ashd{d}$. Fix any generating set $\{{g_1},\ldots,{g_\gamma}\}$ for
$\ash^\ttw/\Rmod$ as an $\ashd{d}$-module. Then $\ash^\ttw/\Rmod$ admits a free $\X$ of rank $\gamma'$ for
some $\gamma'\in\N$ if and only if the generators $\{{g_1},\ldots,{g_\gamma}\}$ admit a matrix of
relations $X(\uld{1})\in\ashd{d}^{\delta\times\gamma}$ which is zero left-prime.
\ET

\BP
Recall the map $\psi:\ashd{d}^\gamma\rightarrow\ash^\ttw/\Rmod$ from equation (\ref{eq:finitely_gen})
and note that $\ash^\ttw/\Rmod$ is isomorphic to $\ashd{d}^\gamma/{\rm ker}(\psi)$ as $\ashd{d}$-modules. Now,
as seen above,
$\ash^\ttw/\Rmod$ admits a free $\X$ of rank $\gamma'$ if and only if $\ash^\ttw/\Rmod$ is a free
$\ashd{d}$-module (of rank $\gamma'$). By
Swan's extension of Serre's conjecture to Laurent polynomial rings, $\ash^\ttw/\Rmod$ is a free
$\ashd{d}$-module if and only if $\ash^\ttw/\Rmod$ is a {\em projective} $\ashd{d}$-module. Since
$\ash^\ttw/\Rmod$ is isomorphic to $\ashd{d}^\gamma/{\rm ker}(\psi)$ as $\ashd{d}$-modules,
$\ash^\ttw/\Rmod$ being projective is equivalent to $\ashd{d}^\gamma/{\rm ker}(\psi)$ being projective.
This, in turn, is equivalent to splitting of the following short exact sequence of $\ashd{d}$-modules
(\cite{eisenbud}):
$$
0\rightarrow {\rm ker}(\psi)\rightarrow \ashd{d}^\gamma\rightarrow\ashd{d}^\gamma/{\rm ker}(\psi)\rightarrow 0.
$$ 
This is true if and only if ${\rm ker}(\psi)$ is a direct summand of $\ashd{d}^\gamma$
(\cite{eisenbud}). But, a submodule of the free module $\ashd{d}^\gamma$ is a direct summand if and only
if it admits representation as the row-span of a zero left-prime matrix (see Remark \ref{rmk:zero:lp}). Thus,
${\rm ker}(\psi)$ being a direct summand of $\ashd{d}^\gamma$ is equivalent to existence of an 
$X(\uld{1})\in\ashd{d}^{\delta\times\gamma}$ such that
$${\rm ker}(\psi)={\rm rowspan}(X(\uld{1})),$$ where $X(\uld{1})$ is zero left-prime.
\EP

\BR
The significance of having a free $\X$ is that initial conditions can then be chosen arbitrarily. That is, solutions
in $\B(\Rmod)$ can then all be obtained by Algorithm \ref{spl-algo} with
arbitrary $x\in\disww{d}{\gamma}$. Note that in Theorem \ref{thm:check:free} it has been assumed that the behavior
is strongly relevant. Therefore, in order to check freeness of initial conditions for a general system, one must first
apply a suitable Noether's normalization, and then invoke Theorem \ref{thm:check:free}. This brings up the question
whether freeness of initial condition is dependent on the particular normalization used, or is it an inherent
property of the original system. This remains an interesting open question. 
\ER

\subsection{Non-autonomy of the set of allowable initial conditions}
As mentioned above, the set of allowable initial conditions, or the state-space, $\X$ can be viewed as
a $d$-D behavior itself. It is often desirable to have this $d$ to be as small as possible. We show in this
subsection that this question of whether $d$ is the smallest possible or not is related with $\X$ being
non-autonomous as a $d$-D behavior.

It is already well-known that an $n$-D behavior $\B$ is non-autonomous if and only if the quotient module
$\ash^\ttw/\Rmod$ is {\em faithful}\footnote{A module $\M$ over $\ash$ is said to be {\em faithful} if ${\rm ann}(\M)
=\{0\}$.} over $\ash$. Translated to the case of $\X$ it means $\X$ is non-autonomous
if and only if $\ash^\ttw/\maphi(\Rmod)$ is faithful as an $\ashd{d}$-module. Keeping this observation in mind,
it is then desirable to design the coordinate transformation $T$ so as to ensure $\ash^\ttw/\maphi(\Rmod)$ is
faithful as an $\ashd{d}$-module. It so happens that our proof of the module theoretic DNNL (Theorem
\ref{thm:mod:dnnl}) already guarantees this property for $\ash^\ttw/\maphi(\Rmod)$.
We bring out this fact in the following Proposition \ref{prop:mod:dnnl}, which is
a slightly modified version of Theorem \ref{thm:mod:dnnl}. 

\BPP\label{prop:mod:dnnl}
Let $\Rmod\subseteq\ash^\ttw$ be a submodule. Suppose the Krull dimension of ${\rm ann}(\ash^\ttw/\Rmod)$ is $d$.
Then there exists a coordinate transformation
$T\in\Z^{n \times n}$ inducing the $\R$-algebra homomorphism $\varphi_T:\ash\ni\uld{}^\nu\mapsto
\uld{}^{T\nu}\in\ash$ for all $\nu\in\Z^n$, and the corresponding $\ash$-module homomorphism
$\maphi:\ash^\ttw\rightarrow\ash^\ttw$,
such that the quotient module $\ash^\ttw/\maphi(\Rmod)$
is a finitely generated and faithful module over $\ashd{d}$.
\EPP

\BP
Recall the proof of Theorem \ref{thm:mod:dnnl}. We have already seen there
that a unimodular $T\in\Z^{n\times n}$ that makes $\ash^\ttw/\maphi(\Rmod)$
a finitely generated module over $\ashd{d}$ can be constructed using Theorem \ref{thm:nnl}: a unimodular
$T\in\Z^{n\times n}$ which renders $\ash/\varphi_T({\rm ann}(\ash^\ttw/\Rmod))$ finitely generated module
over $\ashd{d}$ suffices. However, the $T$ constructed in Theorem \ref{thm:nnl} does more; it makes the canonical
map $\ashd{d}\rightarrow\ash/\varphi_T({\rm ann}(\ash/\Rmod))$ injective and integral. This means the $T$ constructed
in Theorem \ref{thm:nnl} not only makes $\ash/\varphi_T({\rm ann}(\ash/\Rmod))$ finitely generated over $\ashd{d}$,
but faithful over $\ashd{d}$, too. In other words, by Theorem
\ref{thm:nnl}, there exists $T\in\Z^{n\times n}$ unimodular such that $\ash/\varphi_T({\rm ann}(\ash/\Rmod))$
is finitely generated and faithful as a module over $\ashd{d}$. However, $\ash^\ttw/\maphi(\Rmod)$ is clearly
a faithful module over $\ash/\varphi_T({\rm ann}(\ash^\ttw/\Rmod))$. Thus $\ash^\ttw/\maphi(\Rmod)$ is 
faithful over $\ashd{d}$. That $\ash^\ttw/\maphi(\Rmod)$ is finitely generated over $\ashd{d}$ has already been proved
in Theorem \ref{thm:mod:dnnl}.
\EP

Combining this observation with the fact that $\X$ is non-autonomous if and only if $\ash^\ttw/\maphi(\Rmod)$ is
faithful over $\ashd{d}$, we get the following interesting result.

\BC\label{cor:krull:non-auto}
Let $\B\in\Lwe$ be an autonomous $n$-D system such that the Krull dimension of $\ash/{\rm ann}(\ash^\ttw/\Rmod)$
is equal to $d$. Then there exists a unimodular $T\in\Z^{n\times n}$ such that 
\begin{enumerate}
\item $\B$ is $T$-isomorphic to
a strongly relevant autonomous system $\B_2$ of order $d$, and
\item the state space $\X$ corresponding to $\B_2$ is non-autonomous as a $d$-D behavior.
\end{enumerate}
\EC

We have claimed in the beginning of this subsection that non-autonomy of $\X$ is related with the smallest possible $d$.
In order to prove this claim we first show that if the Krull dimension of $\ash/{\rm ann}(\ash^\ttw/\Rmod)$
is equal to $d$ then $\ash^\ttw/\maphi(\Rmod)$ cannot be finitely generated as a module over $\ashd{d'}$ for any
$d'<d$. In other words, for the existence of a $d'$-D behavior $\X$ it is necessary that $d'\geqslant d$.

\BPP\label{prop:krull:non-auto}
Let $\Rmod\subseteq\ash^\ttw$ be a submodule. Suppose the Krull dimension of $\ash/{\rm ann}(\ash^\ttw/\Rmod)$
is equal to $d$. Then there exists a unimodular $T\in\Z^{n\times n}$ such that
$\ash^\ttw/\maphi(\Rmod)$ is finitely generated as a module over $\ashd{d'}$ if and only if $d'\geqslant d$.
\EPP

\BP
{\bf (If)}
We assume that $d'\geqslant d$, and we want to show that there exists a unimodular $T\in\Z^{n\times n}$ such that
$\ash^\ttw/\maphi(\Rmod)$ is finitely generated as a module over $\ashd{d'}$. This is an immediate consequence
of Theorem \ref{thm:mod:dnnl} (or Proposition \ref{prop:mod:dnnl}). Indeed, by Theorem \ref{thm:mod:dnnl} there exists 
a unimodular $T\in\Z^{n\times n}$ such that $\ash^\ttw/\maphi(\Rmod)$ is finitely generated as a module over $\ashd{d}$.
However, since $d'\geqslant d$ means $\ashd{d}$ is a subring of $\ashd{d'}$, $\ash^\ttw/\maphi(\Rmod)$ clearly is
finitely generated as a module over $\ashd{d'}$, too.

{\bf (Only if)}
We assume that there exists a unimodular $T\in\Z^{n\times n}$ such that
$\ash^\ttw/\maphi(\Rmod)$ is finitely generated as a module over $\ashd{d'}$, and we want to show that $d'\geqslant d$.
Let us define $\a:={\rm ann}(\ash^\ttw/\Rmod)$.
Recall from Proposition \ref{prop:ann:id} that $\varphi_T(\a)={\rm ann}(\ash^\ttw/\maphi(\Rmod))$. Since
$\ash^\ttw/\maphi(\Rmod)$ is finitely generated as a module over $\ashd{d'}$, we get that $\ash/{\rm ann}
(\ash^\ttw/\maphi(\Rmod))\a$ is also finitely generated as an $\ashd{d'}$-module. Therefore, by the identity of
Proposition \ref{prop:ann:id}, $\ash/\varphi_T(\a)$ is finitely generated as a module over $\ashd{d'}$. Therefore,
the canonical map $\ashd{d'}/\varphi_T(\a)\cap\ashd{d'}\rightarrow\ash/\varphi_T(\a)$ is injective integral (see
Remark \ref{rmk:integral}). Thus Krull dimension of $\ash/\varphi_T(\a)$ is equal to that of $\ashd{d'}/\varphi_T(\a)
\cap\ashd{d'}$, but Krull dimension of $\ashd{d'}/\varphi_T(\a)\cap\ashd{d'}$ is clearly less than or equal to
that of $\ashd{d'}$,
which is $d'$. Again, Krull dimension of $\ash/\varphi_T(\a)$ is $d$ because $\ash/\varphi_T(\a)$ and $\ash/\a$
are isomorphic (see Remark \ref{rmk:integral}). Therefore, $d'\geqslant d$. 
\EP

Next we show that if $d'>d$ the Krull dimension of $\ash/{\rm ann(\ash^\ttw/\Rmod)}$, then the corresponding
$\X$ cannot be non-autonomous. 

\BL\label{lem:krull:non-auto}
Let $\Rmod\subseteq\ash^\ttw$ be a submodule. Suppose the Krull dimension of $\ash/{\rm ann}(\ash^\ttw/\Rmod)$
is equal to $d$. Further, let $T\in\Z^{n\times n}$ be unimodular such that
$\ash^\ttw/\maphi(\Rmod)$ is finitely generated as a module over $\ashd{d'}$. If $d'>d$ then $\ash^\ttw/\maphi(\Rmod)$
cannot be faithful as a module over $\ashd{d'}$.
\EL

\BP
Suppose, contrary to the claim, that $\ash^\ttw/\maphi(\Rmod)$ is faithful as a module over $\ashd{d'}$. Then
the annihilator of $\ash^\ttw/\maphi(\Rmod)$ as a module over $\ashd{d'}$, which is ${\rm ann}(\ash^\ttw/
\maphi(\Rmod))\cap\ashd{d'}$, must be $\{0\}$. Utilizing the identity of Proposition \ref{prop:ann:id} we get
$$
\varphi_T({\rm ann}(\ash^\ttw/\Rmod))\cap\ashd{d'}=\{0\}.
$$
However, this means that the canonical map $\ashd{d'}\rightarrow\ash/\varphi_T({\rm ann}(\ash^\ttw/\Rmod))$
injective. Further, as has been argued in the proof of Proposition \ref{prop:krull:non-auto}, 
$\ash/\varphi_T({\rm ann}(\ash^\ttw/\Rmod))$ is finitely generated over $\ashd{d'}$ because $\ash^\ttw/\maphi(\Rmod)$
is so. Therefore, the canonical map $\ashd{d'}\rightarrow\ash/\varphi_T({\rm ann}(\ash^\ttw/\Rmod))$ is
integral, too. Therefore, Krull dimension of $\ash/\varphi_T({\rm ann}(\ash^\ttw/\Rmod))$ must be equal to
that of $\ashd{d'}$, which happens to be $d'$. This clearly contradics the assumption that $d'>d$.
\EP

These observations can be summed up as follows:

\BT
Let $\B\in\Lw$ be an autonomous behavior with equation module $\Rmod\subseteq\ash^\ttw$. Let $d$ be the Krull dimension
of $\ash/{\rm ann}(\ash^\ttw/\Rmod)$. Then there exists a unimodular
$T\in\Z^{n\times n}$ such that $\ash^\ttw/\maphi(\Rmod)$ is finitely generated over $\ashd{d'}$ with the corresponding
state space $\X$ non-autonomous if and only if $d'$ is equal to $d$.
\ET

\BP
{\bf (If)} This implication is just a restatement of Corollary \ref{cor:krull:non-auto}.

{\bf (Only if)} By Proposition \ref{prop:krull:non-auto}, $d'\not < d$; for in that case
no unimodular $T\in\Z^{n\times n}$ exists such that $\ash^\ttw/\maphi(\Rmod)$ is finitely generated over $\ashd{d'}$.
On the other hand, if $d'> d$ then, by Lemma \ref{lem:krull:non-auto}, $\ash^\ttw/\maphi(\Rmod)$ cannot be faithful
over $\ashd{d'}$, which means $\X$ cannot be non-autonomous. Therefore, $d'=d$.
\EP

\section{Concluding remarks}
In this paper, we presented novel representation formulae for solutions of discrete $n$-D autonomous systems. Our analysis
is based on a crucial algebraic property: the quotient module being finitely generated over a proper subring of the
ring of shift operators, $\ash$. Quite expectedly, many systems do not have this property. Therefore, we first consider those
systems whose quotient modules are indeed finitely generated over a subring $\ashd{d}$. We called such systems {\em strongly
relevant of order $d$}. We then showed that such systems admit a state space like first order evolutionary equations.
Using this we provided an explicit solution formula for trajectories in such systems. Next, we considered systems that are
not strongly relevant. For these systems we showed that they can be converted to strongly relevant ones by a suitable change
of coordinates on the indexing set $\Z^n$. The crucial fact that enables this conversion is what we have called the discrete
version of Noether's normalization lemma (DNNL). We showed in this paper that given any submodule $\Rmod\subseteq\ash^\ttw$,
there exists a coordinated transformation $T:\Z^n\rightarrow\Z^n$ such that, under the corresponding module homomorphism
$\maphi(\Rmod)$, the quotient module $\ash^\ttw/\maphi(\Rmod)$ becomes finitely generated as a module over a subring $\ashd{d}$.
Using this DNNL we have shown how a general discrete $n$-D autonomous system can be converted to a strongly relevant system.
Then, making use of the earlier representation formula for solutions of strongly relevant systems, and the inverse of the
coordinate transformation we obtained representation formula for general discrete $n$-D autonomous systems.  

An interesting feature of the representation formulae presented in this paper is that the solution trajectories are
obtained by making $(n-d)$-D flow operators acting on initial conditions that are trajectories defined over $\Z^d$. This is
very much akin to the state space solution formula of $1$-D systems. However, there is one crucial deviation: unlike $1$-D
systems, here it may not always possible to chhose the initial conditions freely. In this paper, we provide necessary and sufficient
conditions for these initial conditions to be free. Another interesting question concerning the set of allowable initial conditions
is that whether it can be made smaller. This issue is related with the number $d$. We show in this paper that the smallest value
that the number $d$ can take is equal to the dimension of the characteristic variety of the system.

Many important questions relevant to the results presented in this paper have been left unanswered. For example: can
the representation formulae be extended to non-autonomous systems? Another relevant question is: what is the minimum
number of state variables? Is this number intrinsic to the system, or does it depend on the coordinate transformation?
We have also not answered questions concerning how existing computer algebra packages can be used to compute the various
matrices used in the representation formulae.

\begin{appendix}
\section*{Appendix}
\subsection*{Proof of Lemma \ref{lem:prop:f_g}}
In order to prove Lemma \ref{lem:prop:f_g} we need the following algebraic result: Lemma \ref{lem:app:1}. Lemma
\ref{lem:app:1} follows from elementary commutative algebra. However, for a direct proof of Lemma
\ref{lem:app:1}, the reader may refer to \cite{pal:pillai:sicon}, where $n=2, d=1$ case of the result has been proved.
The extension to the general $n,d$ case is straightforward, and hence we skip it here.

\BL\label{lem:app:1}
The following are equivalent.
\begin{enumerate}
\item $\ash^\ttw/\Rmod$ is finitely generated over $\ashd{d}$. 

\item For every $1\leqslant i \leqslant n-d$ there exists a polynomial $p_i(\uld{})\in{\rm ann}(\ash^\ttw/\Rmod)$
of the form
\begin{equation}
\Dx{d+i}^{L_i}+a_{i,L_i-1}(\uld{1})\Dx{d+i}^{L_i-1}+\cdots+a_{i,1}(\uld{1})\Dx{d+i}+a_{i,0}(\uld{1}),
\end{equation}
where $L_i$ is a finite positive integer, $a_{i,j}(\uld{1})\in\ashd{d}$ for all $1\leqslant i \leqslant n-d$
and $1\leqslant j \leqslant L_i-1$, with $a_{i,0}(\uld{1})$ a unit for all $1\leqslant i \leqslant n-d$.

\item For every $1\leqslant i \leqslant n-d$ there exists a polynomial $p_i(\uld{})\in\ash$ of the form
\begin{equation}
\Dx{d+i}^{L_i}+a_{i,L_i-1}(\uld{1})\Dx{d+i}^{L_i-1}+\cdots+a_{i,1}(\uld{1})\Dx{d+i}+a_{i,0}(\uld{1}),
\end{equation}
where $L_i$ is a finite positive integer, $a_{i,j}(\uld{1})\in\ashd{d}$ for all $1\leqslant i \leqslant n-d$
and $1\leqslant j \leqslant L_i-1$, with $a_{i,0}(\uld{1})$ a unit for all $1\leqslant i \leqslant n-d$, such that
$$
p_i(\uld{})e_j\in\Rmod
$$
for all $1\leqslant j \leqslant \ttw$.
\end{enumerate}
\EL 

With the help of Lemma \ref{lem:app:1} we now prove Lemma \ref{lem:prop:f_g}

\BPref{Lemma \ref{lem:prop:f_g}}
{\bf Part (1)} This assertion follows from repeated application of Lemma \ref{lem:companion}, and the fact that
the maps $\mu_j$ commute pairwise.

{\bf Part (2)} We first prove the result for the special case when $\ttw =1$. 

{\bf \boldmath Case 1 ($q=1$):} In this case $\Rmod\subseteq\ash$ is an
ideal, and hence, ${\rm ann}(\ash/\Rmod)=\Rmod$. Recall that Lemma \ref{lem:prop:f_g} deals with the situation when
$\ash/\Rmod$ is finitely generated as a module over $\ashd{d}$. Therefore, by Lemma \ref{lem:app:1}, for every
$1\leqslant i \leqslant n-d$ there exists a polynomial $p_i(\uld{})\in\Rmod$ having the following form   
$$
\Dx{d+i}^{L_i}+a_{i,L_i-1}(\uld{1})\Dx{d+i}^{L_i-1}+\cdots+a_{i,1}(\uld{1})\Dx{d+i}+a_{i,0}(\uld{1}),
$$
where $L_i$ is a finite positive integer, $a_{i,j}(\uld{1})\in\ashd{d}$ for all $1\leqslant i \leqslant n-d$
and $1\leqslant j \leqslant L_i-1$, with $a_{i,0}(\uld{1})$ a unit in $\ashd{d}$ for all $1\leqslant i \leqslant n-d$.

Now consider a Laurent monomial $\Dx{d+i}^\nu$ where $\nu\in\Z$. Since $p_i(\uld{})$ is monic in $\Dx{d+i}$ and
its coefficient for the zeroth power of $\Dx{d+i}$ is a unit in $\ashd{d}$, we can carry out long division of
$\Dx{d+i}^\nu$ by $p_i(\uld{})$ to get
$$
\Dx{d+i}^\nu=q(\uld{})p_i(\uld{})+b_{i,L_i-1}(\uld{1})\Dx{d+i}^{L_i-1}+b_{i,L_i-2}(\uld{1})\Dx{d+i}^{L_i-2}+\cdots+
b_{i,1}(\uld{1})\Dx{d+i}+b_{i,0}(\uld{1}).
$$
In other words, $\Dx{d+i}^\nu$ is congruent modulo $\Rmod$ to an $\ashd{d}$-linear combination of $\{1, \Dx{d+i},
\Dx{d+i}^2+\cdots+\Dx{d+i}^{L_i-1}\}$. 

\begin{figure}[h!]
\centerline{
\resizebox{7cm}{!}{
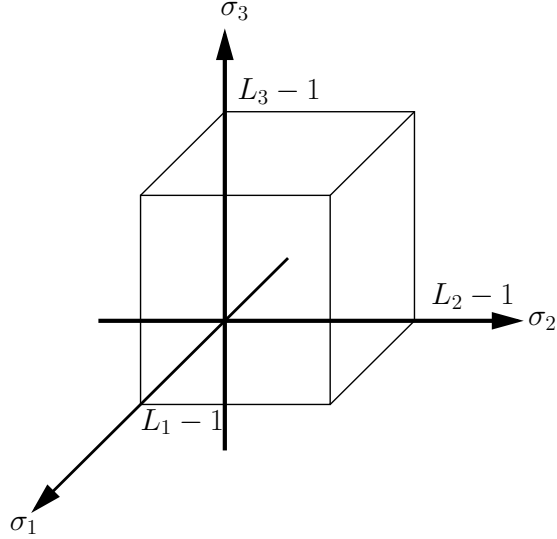
}
}
\caption{The parallelopied when $n-d=3$}\label{fig:pp:grid}
\end{figure}

Notice that the above mentioned procedure can be carried out for all $1\leqslant i \leqslant n-d$. Hence it follows
that a Laurent monomial of the form $\prod_{i=1}^{n-d}\Dx{d+i}^{\nu_i}$ is congruent modulo $\Rmod$ to an
$\ashd{d}$-linear combination of the following monomials:
\begin{equation}
\mathcal{G}:=\left\{\prod_{i=1}^{n-d} \Dx{d+i}^{\nu_i}~\vline~ 0\leqslant \nu_i \leqslant L_i-1\right\}.
\end{equation}
(The indices of the monomials in $\mathcal{G}$ correspond to the parallelopiped, in the positive orthant of the
integer grid $\Z^{n-d}$, with sides having lengths $L_1-1,L_2-1,\ldots, L_{n-d}-1$ units. See Figure
\ref{fig:pp:grid}.) 
Since every Laurent polynomial is an $\ashd{d}$-linear combination of monomials of the form
$\prod_{i=1}^{n-d}\Dx{d+i}^{\nu_i}$, it follows now that every Laurent polynomial is congruent modulo $\Rmod$ to
an $\ashd{d}$-linear combination of monomials in $\mathcal{G}$. 

Hence, the images of the elements in $\mathcal{G}$ under the canonical map $\ash\twoheadrightarrow \ash/\Rmod$
forms a generating set of $\ash/\Rmod$ as a module over $\ashd{d}$. We denote this set of images by $\overline{\G}$.
We claim that under this generating set $\mu_i$'s admit matrix representations that are all unimodular over $\ashd{d}$.

We first show that $\mu_1$ has a matrix representation $A_1(\uld{1})$ that is unimodular. In order to do this we first
put an ordering of the elements in $\G$. Consider the lexicographic ordering of monomials with $\Dx{d+1}\prec\Dx{d+2}
\prec\cdots\prec\Dx{n}$, and write down the elements of $\G$ in descending order. Let us name the elements of
$\overline{\G}$, in this order, as $\{g_1,g_2,\ldots,g_{L}\}$, where $L=\prod_{i=1}^{n-d}L_i$. Now, note that when
the generators are ordered like this then they satisfy the following set of equations:
\begin{equation}\label{eq:app:gen}
\begin{array}{rcl}
\mu_1(g_1) &=& g_2\\
&\vdots&\\
\mu_1(g_{L_{1}}) &=& -a_{1,0}(\uld{1})g_1-a_{1,1}(\uld{1})g_2-\cdots-a_{1,{L_{1}-1}}(\uld{1})g_{L_{1}}\\
\mu_1(g_{L_1+1}) &=& g_{L_1+2}\\
&\vdots&\\
\mu_1(g_{2L_{1}}) &=& -a_{1,0}(\uld{1})g_{L_1}-a_{1,1}(\uld{1})g_{L_1+1}-\cdots-a_{1,{L_{1}-1}}(\uld{1})g_{2L_{1}}\\
&\vdots&
\end{array}
\end{equation} 
It then follows that under this generating set (ordered in the above manner), $\mu_1$ admits a matrix representation
of the following form:
$$
A_1(\uld{1})=\left[\begin{matrix}A_{1,{\rm block}}(\uld{1}) & 0 & \cdots & 0\\
0 & A_{1,{\rm block}}(\uld{1}) & \cdots & 0\\
\vdots & \vdots & \ddots & \vdots\\
0 & 0 & \cdots & A_{1,{\rm block}}(\uld{1}) \end{matrix}\right],
$$
where $A_{1,{\rm block}}(\uld{1})\in\ashd{d}^{L_1\times L_1}$ is given by
$$
A_{1,{\rm block}}(\uld{1})=\left[\begin{matrix}0 & 1 & 0 & \cdots & 0\\
                                               0 & 0 & 1 & \cdots & 0\\
                                               \vdots & \vdots & \vdots & \ddots & \vdots\\
                       -a_{1,0}(\uld{1}) & -a_{1,1}(\uld{1}) & -a_{1,2}(\uld{1}) & \cdots& -a_{1,{L_{1}-1}}(\uld{1})
\end{matrix}\right].
$$
Note that ${\rm det}(A_{1,{\rm block}})=a_{1,0}(\uld{1})$, a unit in $\ashd{d}$. Therefore, $A_{1,{\rm block}}$
is unimodular over $\ashd{d}$. Hence, so is $A_1(\uld{1})$.

We can further claim that not only the map $\mu_1$, in fact, every one of the maps $\mu_2,\mu_3,\ldots,\mu_{n-d}$
has matrix representations under the chosen generating set each of which are unimodular over $\ashd{d}$. In order to
show this let us consider an $i\in\{1,2,\ldots,n-d\}$. Now, we do the following trick: we reorder the elements of
$\G$ so that now once again they are in descending order with a different lexicographic ordering, namely,
$\Dx{d+i}\prec\Dx{d+1}\prec\Dx{d+2}\prec\cdots\prec\Dx{n}$. Once again, we number the elements of $\overline{\G}$
in this order as $\{g'_1,g'_2,\ldots,g'_L\}$ where $L=\prod_{i=1}^{n-d}L_i$. Now, following exactly the same
logic as above we obtain the fact that under this generating set (with this new ordering), the map  $\mu_i$ admits
a matrix representation of the form:
$$
\widetilde{A}_i(\uld{1})=\left[\begin{matrix}A_{i,{\rm block}}(\uld{1}) & 0 & \cdots & 0\\
0 & A_{i,{\rm block}}(\uld{1}) & \cdots & 0\\
\vdots & \vdots & \ddots & \vdots\\
0 & 0 & \cdots & A_{i,{\rm block}}(\uld{1}) \end{matrix}\right],
$$
where $A_{i,{\rm block}}(\uld{1})\in\ashd{d}^{L_i\times L_i}$ is given by
$$
A_{i,{\rm block}}(\uld{1})=\left[\begin{matrix}0 & 1 & 0 & \cdots & 0\\
                                               0 & 0 & 1 & \cdots & 0\\
                                               \vdots & \vdots & \vdots & \ddots & \vdots\\
                       -a_{i,0}(\uld{1}) & -a_{i,1}(\uld{1}) & -a_{i,2}(\uld{1}) & \cdots& -a_{i,{L_{i}-1}}(\uld{1})
\end{matrix}\right].
$$
Therefore, because $a_{i,0}(\uld{1})$ is a unit in $\ashd{d}$, it follows that $A_{i,{\rm block}}(\uld{1})$ is
unimodular over $\ashd{d}$. Hence, $\widetilde{A}_i(\uld{1})$ is also unimodular over $\ashd{d}$. Now note that
the new generating set $\{g'_1,g'_2,\ldots,g'_L\}$ is related with the old generating set $\{g_1,g_2,\ldots,g_L\}$ by a
permutation. Let $P_i$ be the corresponding {\em permutation matrix}. That is,
$$
\left[\begin{matrix}g'_1\\g'_2\\\vdots\\g'_L\end{matrix}\right]=
P_i\left[\begin{matrix}g_1\\g_2\\\vdots\\g_L\end{matrix}\right].
$$
Then, clearly, $P_i^{-1} \widetilde{A}_i(\uld{1})P_i=:
A_i(\uld{1})$ is a matrix representation of the map $\mu_i$ under the old generating set $\{g_1,g_2,\ldots,g_L\}$. The
unimodularity of
$A_i(\uld{1})$ then easily follows from the fact that $\widetilde{A}_i(\uld{1})$ is unimodular.

{\bf \boldmath Case 2 (general $q$):} When $q> 1$ every element in $\ash^\ttw$ can be written as
$$
r_1(\uld{})e_1+r_2(\uld{})e_2+\cdots+r_\ttw(\uld{})e_\ttw,
$$
where $r_i(\uld{})\in\ash$ for all $i$, and $e_j$'s are the standard basis row-vectors of $\ash^\ttw$. Recall that
$\ash^\ttw/\Rmod$ has been assumed to be a finitely generated module over $\ashd{d}$. By Statement 3 of Lemma
\ref{lem:app:1} then
it follows that for every $1\leqslant i < n-d$ there exists a polynomial $p_i(\uld{})\in\ash$ of the form
\begin{equation}\label{eq:app:10}
\Dx{d+i}^{L_i}+a_{i,L_i-1}(\uld{1})\Dx{d+i}^{L_i-1}+\cdots+a_{i,1}(\uld{1})\Dx{d+i}+a_{i,0}(\uld{1}),
\end{equation}
where $L_i$ is a finite positive integer, $a_{i,j}(\uld{1})\in\ashd{d}$ for all $1\leqslant i \leqslant n-d$
and $1\leqslant j \leqslant L_i-1$, with $a_{i,0}(\uld{1})$ a unit for all $1\leqslant i \leqslant n-d$, such that
\begin{equation}\label{eq:app:11}
p_i(\uld{})e_j\in\Rmod
\end{equation}
for all $1\leqslant j \leqslant \ttw$. Suppose we define, like Case 1, the set of monomials 
$$
\G:=\left\{\prod_{i=1}^{n-d}
\Dx{d+i}^{\nu_i}~\vline~ 0\leqslant \nu_i \leqslant L_i-1\right\}.
$$
Then, by equations \eqref{eq:app:10} and
\eqref{eq:app:11}, it follows that for every $1\leqslant i \leqslant \ttw$, $r_i(\uld{})e_i$ is congruent modulo
$\Rmod$ to an $\ashd{d}$-linear combination of elements in $\G$ times $e_i$. In other words, if we define the
following set
$$
\G_{\rm vec}:=\left\{me_j~\vline~ m\in\G, 1\leqslant j \leqslant \ttw\right\},
$$
then the images of $\G_{\rm vec}$ under the canonical map $\ash^\ttw\twoheadrightarrow \ash^\ttw/\Rmod$ is a generating
set for $\ash^\ttw/\Rmod$ as a module over $\ashd{d}$. Now, notice that, under this generating set, the map $\mu_i$
is represented by a matrix
$$
A_{i,{\rm vec}}=\left[\begin{matrix}A_i & 0 & \cdots & 0\\
                                   0 & A_i & \cdots & 0\\
                                   \vdots & \vdots & \ddots & \vdots\\
                                   0 & 0 & \cdots & A_i\end{matrix}\right],
$$
where $A_i$ is the representation of $\mu_i|_{e_j}$. By Case 1, $A_i$ is unimodular. And hence, so is $A_{i,{\rm vec}}$.
 
{\bf Part (3)} Note that if $A_i(\Dx{1})$'s are invertible then the following diagram commutes: 
\begin{equation}\label{eq:commute:2}
\begin{array}{rcl}
\ashd{d}^{\gamma} & \stackrel{\psi}{\twoheadrightarrow} & \M\\
A_i(\uld{1})^{-1}~ \downarrow & & \downarrow~\mu_i^{-1}\\
\ashd{d}^{\gamma} & \stackrel{\psi}{\twoheadrightarrow} & \M
\end{array}.
\end{equation}
Repeated application of Lemma \ref{lem:companion} along with equation \eqref{eq:commute:2} proves that the diagram of
equation \eqref{eq:commute} commutes.

{\bf Part (4)} This assertion follows from Part (3).

\end{appendix}

\nocite{pommaret-quadrat-loc}
\bibliographystyle{alpha}
\bibliography{references}

\begin{thebibliography}{PWR02}

\bibitem[AM69]{atiyah}
M.~F. Atiyah and I.~G. Macdonald.
\newblock {\em Introduction to Commutative Algebra}.
\newblock Addison-{W}esley publishing Co., 1969.

\bibitem[Bj{\"o}79]{bjork}
J.-E. Bj{\"o}rk.
\newblock {\em Ring of Differential Operators}.
\newblock North-Holland Publishing Company, 1979.

\bibitem[Bos82]{bose}
N.K. Bose.
\newblock {\em Applied Multidimensional Systems Theory}.
\newblock Van Nostrand Reinhold, 1982.

\bibitem[CZ95]{curtain}
R.F. Curtain and H.J. Zwart.
\newblock {\em An Introduction to Infinite Dimensional Linear Systems Theory}.
\newblock Springer-Verlag, 1995.

\bibitem[Eis95]{eisenbud}
D.~Eisenbud.
\newblock {\em Commutative Algebra with a View Toward Algebraic Geometry}.
\newblock Springer-Verlag, 1995.

\bibitem[FRZ93]{roc:zam:93}
E.~Fornasini, P.~Rocha, and S.~Zampieri.
\newblock State space realization of 2-{D} finite-dimensional behaviours.
\newblock {\em SIAM Journal on Control and Optimization}, 31(6):1502--1517,
  November 1993.

\bibitem[KR00]{cocoa2}
M.~Kreuzer and L.~Robbiano.
\newblock {\em Computational Commutative Algebra 2}.
\newblock Springer-Verlag, 2000.

\bibitem[NRR11]{avelli:roc:rap:10}
D.~Napp, P.~Rapisarda, and P.~Rocha.
\newblock Time-relevant stability of 2{D} systems.
\newblock {\em Automatica}, 47(11):2373 -- 2382, 2011.

\bibitem[Obe90]{obe:90}
U.~Oberst.
\newblock Multidimensional constant linear systems.
\newblock {\em Acta Appl. Math.}, 20:1--175, 1990.

\bibitem[Par04]{park}
H.~Park.
\newblock Symbolic computation and signal processing.
\newblock {\em Journal of Symbolic Computation}, 37:209--226, 2004.

\bibitem[PP13]{pal:pillai:sicon}
D.~Pal and H.~K. Pillai.
\newblock Representation formulae for discrete 2{D} autonomous systems.
\newblock {\em SIAM Journal on Control and Optimization}, 51(3):2406--2441,
  2013.

\bibitem[PQ99]{pommaret-quadrat-loc}
J-F. Pommaret and A.~Quadrat.
\newblock Localization and parametrization of linear multidimensional control
  systems.
\newblock {\em Systems and Control Letters}, 37:247--260, 1999.

\bibitem[PS98]{hp:shiva:98}
H.~K. Pillai and S.~Shankar.
\newblock A behavioral approach to control of distributed systems.
\newblock {\em SIAM Journal on Control and Optimization}, 37(2):388--408, 1998.

\bibitem[PWR02]{pillai-wood-rogers}
H.~K. Pillai, J.~Wood, and E.~Rogers.
\newblock On homomorphisms of n-{D} behaviors.
\newblock {\em Circuits and Systems I: Fundamental Theory and Applications,
  IEEE Transactions on}, 49(6):732 --742, june 2002.

\bibitem[RW89]{roc:wil:89}
P.~Rocha and J.~C. Willems.
\newblock State for $2$-{D} systems.
\newblock {\em Linear Algebra and its Applications}, 122/123/124:1003--1038,
  1989.

\bibitem[RW97]{state-rapi}
Paolo Rapisarda and J.~C. Willems.
\newblock State maps for linear systems.
\newblock {\em SIAM Journal on Control and Optimization}, 35(3):1053--1091, May
  1997.

\bibitem[Stu02]{sturmfels-poly}
B.~Sturmfels.
\newblock {\em Solving Systems of Polynomial Equations}.
\newblock American Mathematical Society with support from the National Science
  Foundation, 2002.

\bibitem[STW02]{sasane}
A.~J. Sasane, E.~G.~F. Thomas, and J.~C. Willems.
\newblock Time autonomy versus time controllability.
\newblock {\em Systems and Control Letters}, 45:145--153, 2002.

\bibitem[Swa78]{swan}
R.~Swan.
\newblock Projective modules over {L}aurent polynomial rings.
\newblock {\em Transactions of the AMS}, 237:111--120, March 1978.

\bibitem[Val01]{val:01}
M.~E. Valcher.
\newblock Characteristic cones and stability properties of two-dimensional
  autonomous behaviors.
\newblock {\em IEEE Transactions on Circuits and Systems--Part I: Fundamental
  Theory and Applications}, 47(3):290--302, 2001.

\bibitem[Wil83]{will:83}
J.~C. Willems.
\newblock Input-output and state-space representations of finite dimensional
  linear time-invariant systems.
\newblock {\em Linear Algebra and its Applications}, 50:581--608, 1983.

\bibitem[Wil91]{paradigms}
J.~C. Willems.
\newblock Paradigms and puzzles in the theory of dynamical systems.
\newblock {\em IEEE Tansactions on Automatic Control}, 36(3):259--294, March
  1991.

\bibitem[Zam98]{zam:98}
S.~Zampieri.
\newblock Causal input/output representation of 2{D} systems in the behavioral
  approach.
\newblock {\em SIAM Journal on Control and Optimization}, 36(4):1133--1146,
  July 1998.

\bibitem[Zer99]{zerz}
E.~Zerz.
\newblock Primeness of multivariate polynomial matrices.
\newblock {\em Systems and Control Letters}, 29:139--145, 1999.

\end{thebibliography}

\end{document}